\documentclass[review]{elsarticle}               

\usepackage{enumitem}
\usepackage{fixmath}

\usepackage[linesnumbered,norelsize]{algorithm2e}

\newcommand{\newtheoremtmp}[2]{
\newtheorem{#1tmp}{#2}
\newenvironment{#1}{\begin{#1tmp}}{\end{#1tmp}\vspace{1pt}}
}

\usepackage{amsthm}

\newtheoremtmp{thm}{Theorem}
\newtheoremtmp{lem}{Lemma}
\newtheoremtmp{proposition}{Proposition}
\newtheoremtmp{prop}{Property}
\newtheoremtmp{corollary}{Corollary}
\newtheoremtmp{conjecture}{Conjecture}
\newtheoremtmp{condition}{Condition}
\newtheoremtmp{defn}{Definition}
\newtheoremtmp{assum}{Assumption}
\newtheoremtmp{rem}{Remark}
\newtheoremtmp{prob}{Problem}
\newtheoremtmp{exmp}{Example}

\newproof{pf}{Proof}

\usepackage{tikz}
\usetikzlibrary{calc,trees,positioning,arrows,chains,shapes.geometric,%
    decorations.pathreplacing,decorations.pathmorphing,shapes,%
    matrix,shapes.symbols}

\usepackage{pgfplots}
\usepgfplotslibrary{fillbetween}

\pgfplotsset{compat=1.3}

\usepackage{amsmath}
\usepackage{amssymb}
\usepackage{mathtools}

\usepackage{mathrsfs}

\usepackage{caption}
\usepackage{subcaption}

\usepackage{tabularx}

\DeclareMathOperator*{\argmax}{arg\,max}

\newcommand{\T}{\top}

\newcommand{\boundary}{\partial}%

\newcommand\EatDot[1]{}

\newcommand{\Ld}{\mathbf{L}_2}
\newcommand{\Ldloc}{\mathbf{L}_{2,loc}}
\newcommand{\R}{\mathbb{R}}

\newcommand{\Rp}{\R^+}

\newcommand{\SYM}{\mathbb{S}}

\newcommand{\sig}[1]{\mathbold{#1}}
\newcommand{\ssig}[1]{\mathbold{#1}}

\newcommand{\inv}{{\scalebox{0.6}{-1}}}

\newcommand{\sys}{\mathscr{S}}
\newcommand{\Z}{\mathcal{Z}}
\newcommand{\Zp}{\mathcal{Z}_+}
\newcommand{\Zps}{\mathcal{Z}_*}
\newcommand{\N}{\mathbb{N}}

\newcommand{\pp}[1]{\left(#1\right)}
\newcommand{\pb}[1]{\left\{#1\right\}}

\newcommand{\ttime}{\left[0,+\infty\right[}

\newcommand{\ws}{{w^*}}

\newcommand{\lint}{\int\limits}

\newcommand{\smat}[1]{\left[\begin{smallmatrix}#1\end{smallmatrix}\right]}
\newcommand{\bmat}[1]{\begin{bmatrix}#1\end{bmatrix}}

\begin{document}
\begin{frontmatter}
\title{%
Parabolic Set Simulation for Reachability Analysis of Linear Time-Invariant Systems with Integral Quadratic Constraint%
}%
\author[ONERA]{Paul Rousse}\ead{paul.rousse@onera.fr}    
\author[ONERA]{Pierre-Lo\"ic Garoche}\ead{pierre-loic.garoche@onera.fr}               
\author[LAAS-CNRS]{Didier Henrion}\ead{henrion@laas.fr}  

\address[ONERA]{ONERA, 31400 Toulouse, France}  
\address[LAAS-CNRS]{LAAS-CNRS, Universit{\'e} Toulouse, CNRS, Toulouse, France
                and Faculty of Electrical Engineering, Czech Technical University in Prague, Czechia,}             

\begin{abstract} 
This paper describes the computation of reachable sets and tubes for linear
time-invariant systems with an unknown input bounded by integral quadratic
constraints, modeling e.g. delay, rate limiter, or energy bounds. We define a
family of paraboloidal overapproximations. These paraboloids are supported by
the reachable tube on touching trajectories. Parameters of each paraboloid are
expressed as a solution to an initial value problem. Compared to previous
methods based on the classical linear quadratic regulator, our approach can be
applied to unstable systems as well. We tested our approach on large scale
systems.

\end{abstract}
        

\end{frontmatter}

\section{Introduction}
\label{sec:intro}
We consider the reachability problem for Linear Time-Invariant (LTI) with Integral Quadratic Constraints (IQC). Reachable set computation is an active field of research in control theory (see \cite{blanchini2008set}). It has many applications such as state estimation (see \cite{jaulin2001applied}) or verification (see \cite{bayen2007aircraft}) of dynamical systems.
IQC is a classical tool of robust control theory (see e.g. \cite{megretski1997system,megretski2010kyp}). It can model infinite dimensional states, non-linear dynamics, delays, rate limiters, uncertain systems (see \cite{helmersson1999iqc,megretski1997integral,peaucelle2009integral} and \cite{ariba2017}). Up to now, IQCs have mainly been used to evaluate the stability of systems. Despite their ability to model complex systems, we are still lacking results: we do not have a proper characterization of the reachable set.

In this paper, we extend reachability analysis based on ellipsoidal techniques
(see e.g. \cite{chernousko1999, kurzhanski2002ellipsoidal,
kurzhanskiy2007ellipsoidal}) for LTI systems subject to an IQC. This IQC is a
trajectory constraint (i.e. valid at any time) between past state-trajectory,
input signals, and unknown disturbance signals. To override dealing with
constraints over the state-trajectories, we study the LTI system augmented
with a state corresponding to the integral term in the IQC. For a given
parabolic set of initial states, the reachable set of the augmented system is
overapproximated by a time-varying parabolic set. Parameters of this
paraboloid are expressed as the solution to an Initial Value Problem (IVP)
partially described by a Differential Riccati Equation (DRE).  The paraboloid is a
tight overapproximation as it stays in contact with the boundary of the
reachable set on so-called touching trajectories. By studying
touching trajectories that are close to violating the constraint, we find
conditions to generate all the supporting time-varying parabolic sets. At a
given time, the intersection of these supporting parabolic sets is an
exact representation of the reachable set. A wise choice of these
paraboloids can be made to always have a bounded overapproximation of
the reachable set. An algorithm to overapproximate the reachable set is
provided (with an adaptation from \cite{Chandrasekhar} of the Chandrasekhar method for numerical integration
of the DRE in the case where the solution is not
sign-definite). An open source version of our implementation is available on-line \cite{mycode}.

\bigskip
\paragraph*{Related work}
Reachability analysis of LTI systems with ellipsoidal bounded inputs is studied in \cite{chernousko1999,kurzhanski2002ellipsoidal,kurzhanskiy2007ellipsoidal}. Such systems can model infinity norm bounded input-output LTI systems. The reachable set (which is convex and bounded; see \cite{kurzhanski2002ellipsoidal}) can be overapproximated with time-varying ellipsoidal sets. Each ellipsoid is described by its parameters (center and radius) that are solutions to an IVP. These parameters produce tight ellipsoids (i.e., ellipsoids touching the reachable set) which are external approximations of the reachable set. When multiple ellipsoids with different touching trajectories are considered, their intersection is a strictly smaller overapproximation of the reachable set. The accuracy of the overapproximation can be made arbitrarily small by adding more well-chosen ellipsoids. The exact representation of the reachable set is possible by using an uncountable set of ellipsoids.
Our work presents a similar approach for a different class of systems, namely LTI with IQC constraint; and different sets, instead of ellipsoids, we use paraboloids.


An optimal control formulation of the reachable set problem is also possible \cite{lee1967foundations,gusev2017extremal}. For a given state, if the maximal integral cost leading to this state violates the constraint, then this state is unreachable.
It can be solved (using Hamilton-Jacobi-Bellman -HJB- viscosity subsolutions, see \cite{soravia2000viscosity}) leading to global constraints over the reachable set. If the reachable set can be expressed as the intersection (possibly uncountable) of elements of the chosen function family, then the intersection of the resulting constraints gives an exact representation of the reachable set.
However, HJB solutions are known to not scale well with the system dimension.
Our specialized solution showed good results for large systems.

HJB based methods propagate the constraints along the flow of the dynamical system. Occupation measures and barrier certificates methods aim at finding constraints over the reachable tube of a dynamical system: \cite{prajna2004safety} uses IQCs for verification purposes using barrier certificates where the positivity of the energetic state is ensured by using a nonnegative constant multiplier: \cite{henrion2014convex,korda2016moment} use an occupation measure approach where the IQC can potentially be incorporated as a constraint over the moment of the trajectories (note however these references do not deal explicitly with IQCs). A hierarchy of semi-definite programs are derived for polynomial dynamics. Then, off-the-shelf Semi-Definite Programming (SDP) solvers are used to solve the feasibility problem. Optimization-based methods do not usually take advantage of the model structure as they consider a large class of systems (convex, Lipschitz or polynomial dynamics for example).
Similarly than for HJB methods, moment based methods can be used for a large class of systems, but they do not scale well, i.e. they are limited to systems with a small number of states. However, few properties about the reachable set can be formalized such as its domain of existence, its boundedness. We hereby study a narrow class of systems that both have practical and theoretical importance (LTI+IQC system).
Our methods scale favorably with the number of states, and the existence and boundedness of the overapproximation is always granted.

The study of LTI systems with IQC constraint is closely related to the Linear Quadratic Regulator (LQR) problem. In the LQR problem, a quadratic integral is minimized at the terminal time. Optimal trajectories belong to a time-varying parabolic surface, whose quadratic coefficients are a solution to a DRE.
References \cite{savkin1995recursive,guseinov2009approximation,gusev2017extremal} describe the reachable set of LTI systems with terminal IQC.
Reference \cite{jonsson2002robustness} formalizes the problem with a game theory approach. 
Reference \cite{seiler2019finite} solves the differential Riccati inequality over a finite horizon using a basis of polynomial functions, then an SDP solver search for a solution that minimizes the final volume of the overapproximation. This algorithm has been implemented in available tools (see LTVTools toolbox, \cite{LTVTools}).
In all these works, the overapproximation of the reachable set is conditioned by the existence of a solution to the DRA over the interval of integration. In the case of unstable systems, there exists no stable solution to the continuous algebraic Riccati equation. Any reachable set overapproximation is then defined only over a finite interval of time. By taking into account the constraint over the full integration range, we can always find an overapproximation valid over any interval of integration.
Therefore, our method is less conservative than the standard LQR approach.

\paragraph*{Contributions}
We study the reachable set computation of an LTI system with IQC. To the knowledge of the authors, this is the first paper to provide a set-based solution for reachable set computation for LTI systems with IQC that make use of the constraint over the full integration time (not only at the terminal time).
In the conference paper \cite{rousse:hal-02049158}, we presented our preliminary results. Under strong assumptions on the system, we proved that our overapproximation is the exact reachable set computation. In the current paper:
\begin{itemize}
\item We have an exact characterization of the reachable set for a larger class of systems compared to our previous results presented in \cite{rousse:hal-02049158}.
\item We extend the existing ellipsoidal method \citep[presented in][]{chernousko1999,kurzhanski2002ellipsoidal,kurzhanskiy2007ellipsoidal} for the reachability analysis of bounded-input LTI systems to the reachability analysis of LTI systems with IQC. These parabolic constraints are defined by time-varying parameters which are the solution to an IVP. Part of this IVP (the quadratic coefficient of the parabolic constraint) is a DRE. The IVP convergence property is obtained thanks to the convergence property of the DRE.
\item The constraint in the IQC system is modeled as a constraint over the state of an augmented system. Unlike in \cite{CarstenJoost2018,savkin1995recursive,savkin1996model,jonsson2002robustness}, where the constraint is only used at the terminal time, we rather use the constraint on the entire time-domain. The differential equations of the paraboloid's parameters are then differently defined than in previous works. The differential equations depend on a so-called \textit{scaling function}: a time-varying parameter defined by the user. This scaling function can always be chosen such that the overapproximation is defined over any time domain.
\item Our method can efficiently be used for large scale systems for which moment based methods \cite{henrion2014convex,korda2016moment} could not be applied. 
\end{itemize}

\paragraph*{Outline}
The LTI system with IQC and the reachability analysis problem are introduced (Section~\ref{sec:intro}). Parabolic constraints and their associated parameter IVP are defined, their domain of definition is analyzed, the overapproximation property is formulated, as well as the touching trajectories (Section~\ref{sec:paraboloid_over_approx}). A method to generate a set of time-varying parabolic constraints is described. The intersection of these paraboloids exactly describes the reachable set of the system (Section~\ref{sec:exact_reachable_set}). An algorithm to overapproximate the reachable set is described (Section~\ref{sec:implementation}). 
This algorithm is evaluated for numerical examples and large-scale LTI systems (Section~\ref{sec:examples}).

\subsection{Notation}
\newcommand{\mge}{\succ}%
\newcommand{\mgeq}{\succeq}%
\newcommand{\mle}{\prec}%
\newcommand{\mleq}{\preceq}%
\newcommand{\norm}[1]{\left \lVert #1 \right \rVert}
\newcommand{\card}[1]{| #1 |}
\newcommand{\trace}{\mathrm{tr}}
Let $\SYM^{n} \subset \R^{n \times n}$ the set of real valued symmetric square matrices of size $n$. For $A \in \SYM^{n}$, we write $A \mge 0$ (resp. $A \mle 0$) when $A$ is positive definite (resp. negative definite).
We define the matrix norm $\norm{A} = \sqrt{\trace(A^\T A)}$ for $A \in \R^{n \times m}$, where $\trace(B)$ is the trace of $B \in \R^{n \times n}$.
Let a signal be a function that associates to a time instant in $\ttime$ a vector from $\R^n$.
For a given interval $I \subseteq \R$, let $\Ld(I;\R^n)$ denote the Hilbert space of signals equipped with the norm:
\newcommand{\abs}[1]{\left| #1 \right|}
$
\norm{\sig{u}} = \sqrt{\int_{t\in I} \sig{u}^T(t) \sig{u}(t) dt} < \infty.
$
Let $\Ldloc(\Rp;\Rp)$ be the set of locally square integrable signals from $\Rp$ to $\Rp$.
Let $\card{X}$ the cardinality of a countable set $X$.
For a set $\Omega \subset \R^n$, let $\boundary \Omega$ denote its boundary.
\newcommand{\Cont}{\mathscr{C}}%
Let $\Cont^1(I;\R^n)$ the set of functions from $I$ to $\R^n$ which are continuous and differentiable with continuous derivative.
%
\newcommand{\sign}{\mathrm{sign}}
Let $\sign: \R \mapsto \{-1,0,1\}$ such that $\sign(s) = 1$ is $s>0$, $\sign(s) = -1$ if $s<0$, $\sign(s) = 0$ otherwise.


\label{sec:problem_statement}
\newcommand{\extstate}[2]{(#1,#2)}
\newcommand{\Lw}{\Ld(\R^m)}
\newcommand{\Lu}{\Ld(\R^p)}
\newcommand{\Lx}{\Ld(\R^m)}
\subsection{System}
For a given input signal $\sig{u} \in \Cont^1(\Rp;\R^p) \cap \Ldloc(\Rp;\R^p)$, given matrices $A\in\R^{n \times n}$, $B \in \R^{n \times m}$, $B_u \in \R^{n \times p}$, and a given terminal time $t>0$, we study the trajectories $\sig{x} \in \Ld([0,t];\R^n)$ of the LTI system:
\newcommand{\td}{\tau}
\newcommand{\tint}{s}
\begin{equation}
\label{eq:ode}
\left\{
\begin{aligned}
\dot{\sig{x}}(\td) &= A \sig{x}(\td) + B \sig{w}(\td) + B_u \sig{u}(\td) && \textrm{with $\td \in [0,t]$}\\
\sig{x}(0) &= x_0
\end{aligned}
\right.
\end{equation}
where $\sig{w} \in \Ldloc([0,t];\R^m)$ is an unknown disturbance that satisfies:
\begin{equation}
\label{eq:cons}
x_{q0} +
\lint_0^\td
\left[\begin{smallmatrix}\sig{x}(\tint)\\\sig{u}(\tint)\\\sig{w}(\tint)\end{smallmatrix}\right]^\top
M
\left[\begin{smallmatrix}\sig{x}(\tint)\\\sig{u}(\tint)\\\sig{w}(\tint)\end{smallmatrix}\right]
d\tint \geq 0 \textrm{ for all $\td \in [0,t]$}
\end{equation}
for given initial conditions $\extstate{x_0}{x_{q0}} \in \R^n \times \Rp$, and given symmetric matrix
\begin{equation}
\label{def:M_matrix}
M = 
\begin{bmatrix}
M_x & M_{xu} & M_{xw} \\ M_{xu}^\T & M_{u} & M_{uw} \\ M_{xw}^\T &  M_{uw}^\T & M_{w} 
\end{bmatrix} \in \SYM^{n+m+p}
\end{equation}
with $M_{w} \mle 0$.
Many IQC models verify such condition, as the one described in \cite{megretski1997system}.

In this work, the constraint \eqref{eq:cons} is expressed as a constraint over a state $\sig{x_q} \in \Ld([0,t];\R)$ defined for $s \in [0,t]$ by:
\begin{equation}
\label{eq:xq_int}
\sig{x_q}(\td) = 
x_{q0} +
\lint_0^\td
\left[\begin{smallmatrix}\sig{x}(\tint)\\\sig{u}(\tint)\\\sig{w}(\tint)\end{smallmatrix}\right]^\top
M
\left[\begin{smallmatrix}\sig{x}(\tint)\\\sig{u}(\tint)\\\sig{w}(\tint)\end{smallmatrix}\right]
d\tint,
\end{equation}
then
\begin{equation}
\label{eq:state_cons}
\sig{x_q}(\td) \geq 0 \textrm{ for all $\td \in [0,t]$.}
\end{equation} 
The constrained dynamical system $\sys(\Z_0,t)$ is then defined for a given set of initial states $\Z_0 \subset \R^n \times \R$ and a terminal time $t>0$:
\begin{equation}
\label{def:cons_sys}
\sig{z} = \extstate{\sig{x}}{\sig{x_q}} \in \sys(\Z_0,t)
\Leftrightarrow
\left\{
\begin{array}{l}
\textrm{$\sig{x}$ solves \eqref{eq:ode}}\\
\textrm{and $\sig{x_q}$ solves \eqref{eq:xq_int}}\\
\textrm{with $\extstate{x_0}{x_{q0}} \in \Z_0$}\\
\textrm{$\sig{x_q}$ satisfies \eqref{eq:state_cons}}
\end{array}
\right.
\end{equation}
\newcommand{\reach}{\mathcal{R}}%
Define the reachable set:
\begin{equation}
\label{eq:reach_set}
\reach(\Z_0,t) = \pb{\sig{z}(t) \middle| \sig{z} \in \sys(\Z_0,t)}.
\end{equation}
Then, $\reach(\Z_0,t) \subseteq \Zp$ where $\Zp = \R^n \times \Rp$,
and let $ \Zps = \R^n \times \{0\}$.

\newcommand{\parset}{\Pi}%
\newcommand{\parabol}{\mathcal{P}}%
\newcommand{\parinter}{\parset^\cap}%
\newcommand{\parabolSet}{\mathbb{P}}%
\newcommand{\fc}{\ssig{f}}
\newcommand{\fct}{\fc^\T}
\renewcommand{\r}{\ssig{g}}
\newcommand{\E}{\ssig{E}}

\newcommand{\fci}{f}
\newcommand{\ri}{g}
\newcommand{\Ei}{E}

\newcommand{\dfc}{\dot{\fc}}
\newcommand{\dr}{\dot{\r}}
\newcommand{\dE}{\dot{\E}}

\subsection{Paraboloids} 
We overapproximate the reachable set $\reach(\Z_0,t)$ of $\sys(\Z_0,t)$ with \textit{paraboloids}:
\begin{defn}[Paraboloid]
\label{def:semiparabol}
Given $(\Ei,\fci,\ri) \in \SYM^{n} \times \R^n \times \R$, define the \textit{value function}:
\begin{equation*}
\begin{array}{rcl}
h : &\R^n \times \R &\to \R\\
    &\extstate{x}{x_q} & \mapsto  x^\T \Ei x - 2 \fci^\T x + \ri +x_q,
\end{array}
\end{equation*}
and the paraboloid:%
\begin{equation*}
\label{eq:semiparabol}
\parabol(\Ei,\fci,\ri) = \pb{ \extstate{x}{x_q} \in \R^{n+1} \middle| h\extstate{x}{x_q} \leq 0}.
\end{equation*}
\end{defn}

\newcommand{\scaling}{\gamma}%
\newcommand{\dxqs}{\sig{\dxq^*}}%
\newcommand{\ddxqs}{\ddot{x}_q^*}%
\newcommand{\tp}{l}%
Let $\parabolSet = \pb{\parabol(\Ei,\fci,\ri) \middle | \Ei \in \SYM^n, \fci \in \R^n, \ri \in \R}$ be the set of paraboloids. Elements of $\parabolSet$ are not strictly speaking paraboloids since $E \mgeq 0$ is not assumed in the definition.
\begin{defn}[Scaled Paraboloid]
\label{def:scaled_paraboloid}
For $\parabol \in \parabolSet$ with parameters $(\Ei,\fci,\ri)$ and a scaling factor $\scaling > 0$, let $\scaling \parabol \in \parabolSet$ be the scaled paraboloid defined by parameters $(\scaling \Ei,\scaling \fci,\scaling \ri)$.
\end{defn}%
Scaled paraboloids satisfy the following:
\begin{prop}
\label{prop:scaled_overapprox}
Given $\parabol \in \parabolSet$ and $\scaling \geq 1$, it holds $\parabol \cap \Zp \subseteq \scaling \parabol \cap \Zp$.
\end{prop}
\begin{pf}
Let $h$ and $h'$ (resp.) the value functions of $(\Ei,\fci,\ri) = \parabol$ and $\scaling \parabol$ (resp.) evaluated at $\extstate{x}{x_q} \in \parabol$. Since $\extstate{x}{x_q} \in \parabol$, $h \leq 0$, i.e. $x^\T \Ei x - 2 \fci^\T x + \ri  \leq -x_q$. Then, $h' =  \scaling (x^\T \Ei x - 2 \fci^\T x + \ri) + x_q  \leq - (\scaling-1) x_q$. Since $\extstate{x}{x_q} \in \Zp$ and since $\scaling-1 \geq 0$, we have $(\scaling-1) x_q \geq 0$ i.e. $h' \leq 0$ meaning that $\extstate{x}{x_q} \in \scaling \parabol \cap \Zp$.
\end{pf}

For $P$ a function that associates to a time $t$ of a time-interval $I \subset \Rp$ a set of states $P(t) \subset \R^{n+1}$. We define a \textit{touching trajectory}:
\begin{defn}[Touching Trajectory]
A trajectory $\sig{z}^*$ solution to (\ref{eq:ode}, \ref{eq:xq_int}) is a touching trajectory of $P$ when $\sig{z}^*(t)$ belongs to the surface of $P(t)$ at any time $t \in I$, i.e. $\sig{z}^*(t) \in \boundary P(t)$.
\end{defn}%

\subsection{Problem Statement}
We are now ready to state the two problems studied in this work.
\begin{prob}
\label{pb:approx_reach_set}
Find an overapproximation of the reachable set $\reach(\parabol_0,t)$ at any $t>0$ for a given paraboloid of initial conditions $\parabol_0 \in \parabolSet$.
\end{prob}
Theorem~\ref{thm:overapproximation} in Section~\ref{sec:paraboloid_over_approx} solves Problem~\ref{pb:approx_reach_set}. It restates classical results about LQR systems applied to reachability analysis of IQC systems when the constraint \eqref{eq:cons} is a terminal time constraint. Theorem~\ref{thm:overapprox_scaled} in Section~\ref{sec:exact_reachable_set} is another solution to Problem~\ref{pb:approx_reach_set}. It extends the result in Theorem~\ref{thm:overapproximation} taking into account the constraint over the entire interval of integration.

\begin{prob}
\label{pb:exact_reach_set}
Find a sequence of overapproximations that converges to the reachable set $\reach(\parabol_0,t)$ at any $t>0$ for a given paraboloid of initial conditions  $\parabol_0 \in \parabolSet$.
\end{prob}
Theorem~\ref{thm:exact_reachable_set} in Section~\ref{sec:exact_reachable_set} solves Problem~\ref{pb:exact_reach_set}.
We prove for any given state on the boundary of $\reach(\parabol_0,t)$ that there exists a tight overapproximation touching $\reach(\parabol_0,t)$ on this state.
Then, we can describe the reachable set as the intersection of all the possible overapproximations.

\section{Overapproximation with Paraboloids}
\label{sec:paraboloid_over_approx}
In this section, Problem~\ref{pb:approx_reach_set} is solved using time-varying paraboloids $P: I\rightarrow \parabolSet$ where $I$ is the interval of definition of $P$. Time-varying parameters $(\E,\fc,\r)$ of $P$ satisfies a differential equation that guarantees an overapproximation relationship with the reachable set, i.e. $\reach(\parabol_0,t) \subseteq P(t)$ for any $t \in I$. We express existence and domain of definition $I$ of the time-varying paraboloid $P$. We prove that the overapproximations $P$ are \textit{tight} since there are touching trajectories of $\reach(\parabol_0,t)$ that both belong to the surface of $P(t)$ and to the surface of $\reach(\parabol_0,t)$ for $t\in I$. Finally, the method is presented for a simple toy example.

\bigskip

\newcommand{\qi}{q^\inv}%
\newcommand{\Msc}{M^{sc}}%
\newcommand{\Mw}{M_{w}}%
\newcommand{\Mwi}{\Mw^\inv}%
\newcommand{\Mxw}{M_{xw}}%
\newcommand{\Mx}{M_{x}}%
\newcommand{\Mxwt}{M_{xw}^\T}%
\newcommand{\Qi}{Q^\inv}%
\newcommand{\xw}{\smat{x\\w}}%
\newcommand{\xwt}{\smat{x\\w}^\T}%
\newcommand{\Bt}{B^\T}%
\newcommand{\At}{A^\T}%
\newcommand{\xc}{x_c}%
\newcommand{\xct}{x_c^\T}%
\newcommand{\dxc}{\dot{x}_c}%
\newcommand{\xq}{\sig{x_q}}%
\newcommand{\dxq}{\dot{x}_q}%
\newcommand{\dq}{\dot{q}}%
\newcommand{\dQ}{\dot{Q}}%
\newcommand{\xxc}{\pp{x-\xc}}%
\newcommand{\xxct}{\xxc^\T}%
\newcommand{\Px}{p_x}%
\newcommand{\Pu}{p_u}%
\newcommand{\PxMsc}{\Px \Msc}%
\newcommand{\PxMscPx}{\Px^\T \Msc \Px}%
\newcommand{\xu}{\smat{x\\u}}%
\newcommand{\xut}{\smat{x\\u}^\T}%
\newcommand{\xcu}{\smat{\xc\\u}}%
\newcommand{\xcut}{\smat{\xc\\u}^\T}%
\newcommand{\xuw}{\smat{x\\u\\w}}%
\newcommand{\xuwt}{\smat{x\\u\\w}^\T}%
\newcommand{\Bu}{B_u}%
\newcommand{\Muw}{M_{uw}}%
\newcommand{\Muwt}{\Muw^\T}%
\newcommand{\Mu}{M_{u}}%
\newcommand{\Mxu}{M_{xu}}%
\newcommand{\Mscu}{M^{sc}_u}%
\newcommand{\pu}{p_u}%
\newcommand{\Mxuwt}{\smat{ \Mxw \\ \Muw}^\T}%
\newcommand{\xuws}{\smat{x\\u\\\ws}}%
\newcommand{\xuwst}{\xuws^\T}%
%
%
Parameters of $P$ are expressed as solutions to an initial value problem. For given $\Ei_0 \in \SYM^n$, let $\E$ be the solution to the following Differential Riccati Equation (DRE) with initial condition $\E(0) = \Ei_0$:
\begin{equation}
\label{eq:eq_diff_E}
\scalebox{0.99}{\mbox{\ensuremath{\displaystyle
\begin{aligned}
\dE(t) = & - \E(t) A - A^\T \E(t) - \Mx\\
& + \pp{\Bt \E(t) + \Mxwt}^\T \Mwi \pp{\Bt \E(t) + \Mxwt}.
\end{aligned}
}}}
\end{equation}%
\newcommand{\Tdef}[1]{T_{P}(#1)}%
\newcommand{\TdefE}[1]{T_{E}(#1)}%
\newcommand{\IdefE}[1]{\left[0,\TdefE{#1}\right[}%
\newcommand{\Idef}[1]{\mathcal{I}(#1)}%
Let $\TdefE{\Ei_0} \in \Rp \cup \{+\infty\}$ be  defined for the initial condition $\Ei_0$ s.t. $\IdefE{\Ei_0}$ is the interval of definition of the solution of \eqref{eq:eq_diff_E} (existence, uniqueness, convergence properties and continuity of the solution are studied in \cite{kuvcera1973review}).
Let $\fc$ denote the solution to the following IVP with initial condition $\fc(0) = \fci_0$:
\begin{equation}
\label{eq:eq_diff_fc}
\begin{split}
\dfc(t) = &- A^\T \fc(t) + (\Mxu + \E(t) \Bu )\sig{u}(t)\\
&+ (\E(t) B + \Mxw ) \Mwi (\Bt \fc(t) -\Muwt \sig{u}(t) ).
\end{split}
\end{equation}
$\fc$ satisfies a linear varying parametric differential equation with a continuous input signal. On $\IdefE{\Ei_0}$, solution $\fc$ to \eqref{eq:eq_diff_fc} exists, is unique and continuous.
By continuity of $\fc$ and $\sig{u}$ over $\IdefE{\Ei_0}$, $\r$ is defined on $\IdefE{\Ei_0}$.
For $t \in \IdefE{\Ei_0}$, let:
\begin{equation}
\label{eq:eq_diff_g}
\r(t) = \ri_0 + \lint_{0}^t \smat{\fc(\tau)\\\sig{u}(\tau)}^\T G \smat{\fc(\tau)\\\sig{u}(\tau)} d\tau
\end{equation}
where
$$
G = \begin{bmatrix}
B \Mwi \Bt & \Bu - B \Mwi \Muwt\\ (\Bu - B \Mwi \Muwt)^\T & -\Mu + \Muw \Mwi \Muwt
\end{bmatrix}.
$$

\newcommand{\feqdiff}{F}
\newcommand{\fE}{\mathcal{F}_{\E}}
\newcommand{\ffcA}{\mathcal{A}_{\fc}}
\newcommand{\ffcB}{\mathcal{B}_{\fc}}
\newcommand{\fr}{\mathcal{F}_{\r}}
\newcommand{\IVP}{\mathcal{T}}
\begin{defn}[Time-Varying Paraboloid]
\label{def:eq_diff_well_defined}
For an initial paraboloid $\parabol_0 \in \parabolSet$, let the time-varying paraboloid $P$ be defined as:
\[
\begin{array}{rl}
P \colon I &\to \parabolSet\\
   t &\mapsto \parabol(\E(t),\fc(t),\r(t))
\end{array}
\]
where the time-varying coefficients $(\E,\fc,\r)$ are solutions to (\ref{eq:eq_diff_E},\ref{eq:eq_diff_fc},\ref{eq:eq_diff_g}) with initial condition $\parabol(\Ei_0,\fci_0,\ri_0) = \parabol_0
$.
Let $\IVP$ be the function that associates to the initial paraboloid $\parabol_0 \in \parabolSet$ the time-varying paraboloid $P$. Let $\Tdef{P} = \TdefE{\Ei_0}$ and $\Idef{P} = \left[0, \Tdef{P}\right[$ be the interval of definition of $P$.
\end{defn}
%

\newcommand{\htX}{h_{z}}
\newcommand{\htXs}{h_{z^*}}
\newcommand{\dhtX}{\dot{h}_{z}}
\newcommand{\dhtXs}{\dot{h}_{z^*}}
For $P = \IVP(\parabol_0)$, let $h(t,\cdot)$ be the value function of $P(t)$ at $t \in \Idef{P}$. For $z_t = \extstate{x_t}{x_{q,t}} \in \R^{n+1}$, $w_t \in \R^m$, $\htX(t) = h(t,\sig{z}(t))$ is the value function along the trajectory $\sig{z} = \extstate{\sig{x}}{\sig{x_q}}$ solution to (\ref{eq:ode}, \ref{eq:xq_int}) generated by $\sig{w}$ such that $\sig{w}(t) = w_t$ and $\sig{z}(t) = z_t$.
\begin{prop}
\label{prop:max_dhX}
The maximum time derivative of the value function $h(t,\sig{z}(t))$ along the trajectories $\sig{z}$ for a disturbance $w_t$ exists for all $t \in \Idef{P}$ and it is equal to zero.
\end{prop}
\begin{pf}
At a time $t \in \Idef{P}$, it holds
$$
\htX(t) = 
\bmat{\sig{x}(t)\\1}^\T 
\bmat{\E(t) & -\fc(t)\\-\fc^\T(t) & \r(t)}
\bmat{\sig{x}(t)\\1}^\T
+ \xq(t)
$$
and the time derivative of $\htX$ at $t$ is
\begin{equation}
\label{eq:dhtX}
\begin{aligned}
\dhtX(t) &= 
\smat{\sig{x}(t)\\1}^\T
\smat{\dE(t) & -\dfc(t)\\-\dfc^\T(t) & \dr(t)}
\smat{\sig{x}(t)\\1}
\\
&+
2\smat{\sig{x}(t)\\1}^\T
\smat{\E(t) & -\fc(t)\\-\fc^\T(t) & \r(t)}
\smat{\dot{\sig{x}}(t)\\0}
\\&
+ 
\smat{\sig{x}(t)\\\sig{u}(t) \\ w_t}^\T M 
\smat{\sig{x}(t)\\\sig{u}(t) \\ w_t}.
\end{aligned}
\end{equation}
Therefore, $k(w_t) = \dhtX(t)$ is a quadratic function of $w_t$:
$$
\begin{aligned}
k(w_s) =& \sig{d}(t)^\T + 2 \pp{\Bt (\E \sig{x}-\fc) + \Mxuwt \smat{\sig{x}\\\sig{u}}}^\T w_t 
\\&+ w_t^\T M_w w_t
\end{aligned}
$$
where $\sig{d}$ is function of $\E(t)$, $\fc(t)$, $\r(t)$, $\sig{x}(t)$ and $\sig{u}(t)$ and system parameters. Since $\Mw \mle 0$, the supremum of $w_t \mapsto k(w_t)$ exists and is attained for $w_t = \sig{\ws}(t) = \argmax_{w_t \in \R^m}  k(w_t)$ with:
\begin{equation}
\label{def:ws}
\sig{\ws} = - \Mwi \pp{\Bt (\E \sig{x}-\fc) + \Mxuwt \smat{\sig{x}\\\sig{u}}}.
\end{equation}
Since $M_w \mle 0$, $\Mwi$ is well defined. Using (\ref{eq:eq_diff_E},\ref{eq:eq_diff_fc},\ref{eq:eq_diff_g}) in \eqref{eq:dhtX}, we get $\max_{w_t \in \R^m}  \dhtX = 0$. Therefore, $\dhtX \leq 0$ for any $w_t \in \R^m$. 
\end{pf}

\bigskip
We can now state one of our main results:
\begin{thm}[Solution to Problem~\ref{pb:approx_reach_set}]
\label{thm:overapproximation} Let $P = \IVP(\parabol_0)$ for a set of initial states $\parabol_0$. For all $t \in \Idef{P}$, the reachable set $\reach(\parabol_0,t)$ of $\sys(\parabol_0,t)$, is overapproximated by $P(t)$, i.e.:
\[
\forall t \in \Idef{P}, \reach(\parabol_0,t) \subseteq P(t) \cap \Zp
.\]
\end{thm}%
\begin{pf}
\label{subsec:eq_diff_uncentered}
Using Property~\ref{prop:max_dhX}, by integration of $\dhtX$, if $\htX(0) \leq 0$ then $ \forall t \in \Idef{P}, \htX(t) \leq 0$, i.e.:
$$\sig{z}(0) \in P(0)  \Rightarrow \sig{z}(t) \in P(t) \textrm{ for all $t \in \Idef{P}$}.$$
The constraint \eqref{eq:state_cons} ensures that $\sig{z}(t) \in \Zp$.
\end{pf}

\begin{prop}
Let $\sig{z^*}$ be a trajectory generated by $\sig{\ws}$ defined in \eqref{def:ws} such that initial condition satisfies $\sig{z^*}(0) \in \boundary \parabol_0$. At any time $t\in \Idef{P}$, it holds  $\sig{z^*}(t) \in \boundary P(t)$.
\end{prop}%
\begin{pf}
$\sig{z^*}$ is the trajectory generated by the optimal disturbance $\sig{\ws}$. Using Property~\ref{prop:max_dhX}, $\dhtXs(t) = 0$ for any $t\geq 0$. Since $\htXs(0) = 0$, by integration, $\htXs(t) = 0$.
\end{pf}
Trajectories generated by $\sig{\ws}$ defined in \eqref{def:ws} stay in contact with the surface of their time-varying paraboloids. Touching trajectories of $P$ do not necessarily belong to $\sys(\parabol_0,t)$, $t \in \Idef{P}$, as the energetic constraint might be locally violated.

\begin{rem}[Representation of paraboloids]
In \cite{savkin1995recursive}, the time-varying value function is a quadratic function defined by its quadratic coefficient $\ssig{S}$, its center $\ssig{x_c}$ and its value at the center $\ssig{\rho}$, all satisfying an IVP. In \cite{savkin1995recursive}, the center $\ssig{x_c}$ can diverge when the determinant of $\ssig{S}$ vanishes. However, the corresponding time-varying value function is time-continuous and can be extended continuously. In this paper, we choose to work with variables $\E$, $\fc$ and $\r$ (see Definition~\ref{def:semiparabol}) to avoid this issue.
\end{rem}

\begin{figure*}
\newcommand{\imgLarge}{0.23\textwidth}
\centering
\includegraphics[angle=90,scale=0.4]{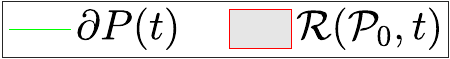}
\includegraphics[width=\imgLarge]{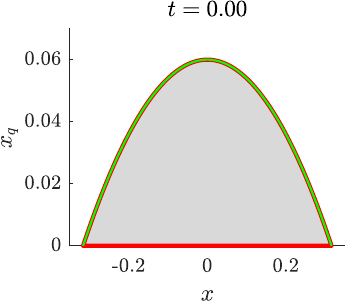}
\includegraphics[width=\imgLarge]{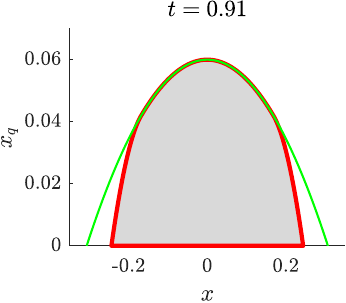}
\includegraphics[width=\imgLarge]{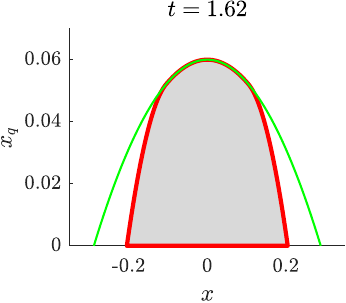}
\includegraphics[width=\imgLarge]{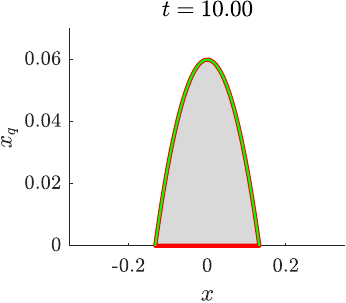}
\caption{\label{fig:bounding_paraboloid_high_energy}%
Time-varying paraboloid overapproximating the reachable set at different time instants $t$ in $\{0.00,0.91,1.62,10.00\}$ for an initial maximum energetic level of $x_{q,0} = 0.06$. The solution to \eqref{eq:eq_diff_E} converges to a constant value when $t \rightarrow +\infty$. The shaded regions are the reachable set $\reach(\parabol_0,t)$, the thin lines are the boundary of the overapproximation $P(t)$ of Theorem~\ref{thm:overapproximation}.}
\end{figure*}

\begin{figure*}
\newcommand{\imgLarge}{0.23\textwidth}
\centering
\includegraphics[angle=90,scale=0.4]{img/dummy_legend.pdf}
\includegraphics[width=\imgLarge]{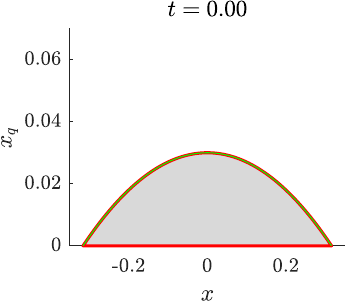}
\includegraphics[width=\imgLarge]{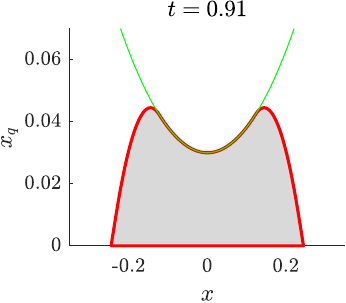}
\includegraphics[width=\imgLarge]{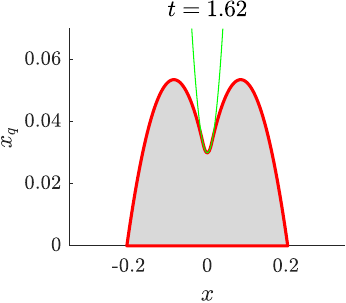}
\includegraphics[width=\imgLarge]{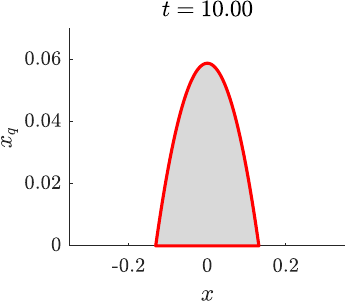}
\caption{\label{fig:bounding_paraboloid_low_energy}%
Time-varying paraboloid overapproximating the reachable set at different time instants $t$ in $\{0.00,0.91,1.62,10.00\}$ for an initial maximum energetic level of $x_{q,0} = 0.03$. The solution to \eqref{eq:eq_diff_E} has a finite escape time and diverge at $t = 1.68$. The shaded regions are the reachable set $\reach(\parabol_0,t)$, the thin lines are the boundary of the overapproximation $P(t)$ of Theorem~\ref{thm:overapproximation}.}
\end{figure*}

\begin{exmp}
\label{ex:1dsys}
Let
$A = -1$, 
$B = 1$,
$M = \smat{1&0&0\\0&1&0\\0&0&-2}$,
$B_u = 0$ and 
$\sig{u}: \left[0,\infty\right[ \mapsto 0$.
Solutions to IVP~\eqref{eq:eq_diff_E} (that is $\dE =-\frac{1}{2} \E^2 + 2 \E -1$) diverge for $\Ei_0 \mle \Ei^-$ (see Figure~\ref{fig:stab_care}) where $\Ei^- \mle \Ei^+$ are the roots of the equation $-\frac{1}{2} \Ei^2 + 2 \Ei -1 = 0$ for $\Ei \in \R$, $\Ei^- = 2-\sqrt{2}$ and $\Ei^+ = 2+\sqrt{2}$.
\begin{figure}
\centering
\includegraphics[width=0.8\columnwidth]{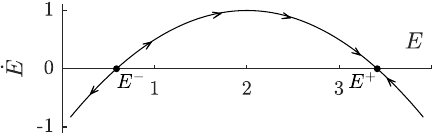}
\caption{\label{fig:stab_care}Convergence analysis of the DRE for Example~\ref{ex:1dsys}}
\end{figure}
Figure~\ref{fig:bounding_paraboloid_high_energy} shows the trajectory of the paraboloid for $\Ei_0$ in the stable region $\Ei_0 \mge \Ei^-$ while Figure~\ref{fig:bounding_paraboloid_low_energy} shows the trajectory of the paraboloid for $\Ei_0$ in the unstable region $\Ei_0 \mle \Ei^-$.
\end{exmp}

\subsection{Domain of definition of $P$}
By Definition~\ref{def:eq_diff_well_defined}, the domain of a time-varying paraboloid $P$ is the domain of its quadratic time-varying coefficient $\E$.
Since the solution of the DRE~\eqref{eq:eq_diff_E} might diverge in a finite time $\TdefE{\Ei_0} < \infty$ (where $\Ei_0$ is the initial condition of \eqref{eq:eq_diff_E}), $P$ is defined only in the right-open interval $\left[0,\TdefE{\Ei_0}\right[$. In this part, we show that since the touching trajectories of $P$ are  defined over the closed interval $\left[0,\TdefE{\Ei_0}\right]$,the definition of $P$ can be prolongated to the same closed interval.

Property~\ref{prop:max_dhX} and (\ref{eq:eq_diff_E},\ref{eq:eq_diff_fc},\ref{eq:eq_diff_g}) can be derived solving the following optimal control problem (for $t>0$):
\begin{equation*}
\begin{array}{rl}
\displaystyle\max_{\sig{w} \in \Ld([0,t];\R^m)}{} & \displaystyle\lint_0^t \smat{\sig{x}(\tau)\\\sig{u}(\tau)\\\sig{w}(\tau)} M \smat{\sig{x}(\tau)\\\sig{u}(\tau)\\\sig{w}(\tau)} d\tau - x_{q,t}\\
\textrm{s.t.} & \dot{\sig{x}} = A \sig{x} + B \sig{w} + B_u \sig{u}\\
 & \sig{x}(t) = x_t\\
\end{array}
\end{equation*}
for given $\extstate{x_t}{x_{q,t}} \in \Zp$. This is a special instance of the LQR problem (see e.g. \cite{savkin1995recursive}).
For $x \in \Ld(T;\R^n)$ a touching trajectory,  let
$$
\sig{n} = \E \sig{x} - \fc.
$$
be the normal to the paraboloid surface.
Using (\ref{eq:eq_diff_E},\ref{eq:eq_diff_fc},\ref{eq:eq_diff_g}), $\sig{n}$ satisfies the following differential equation:
$$
\begin{bmatrix}
\dot{\sig{x}}\\ \dot{\sig{n}} \end{bmatrix}
 = 
L
\begin{bmatrix}
\sig{x} \\ \sig{n}
\end{bmatrix}
+
N\sig{u}
$$%
\newcommand{\Ah}{L}%
\newcommand{\Bh}{N}%
where
$$
\Ah
=
\begin{bmatrix}
A - B \Mwi \Mxwt & - B \Mwi \Bt \\
-(\Mx -  \Mxw \Mwi \Mxwt) & -\At + \Mxw \Mwi \Bt
\end{bmatrix}
$$
and
$$\Bh
=
\begin{bmatrix}
\Bu - B \Mwi \Muwt\\
-(\Mxu - \Mxw \Mwi \Muwt)
\end{bmatrix}.
$$

The value function evaluated along the touching trajectory $\sig{x}$ is then obtained by introducing the parameter $\sig{r} = \r - \fct \sig{x}$ which satisfies:
$$
\dot{\sig{r}} = \sig{u} \begin{pmatrix} H & R \end{pmatrix} \smat{\sig{x} \\ \sig{n} \\ \sig{u}}
$$
with
$$
H = \begin{pmatrix}\Muw \Mwi \Mxw - \Mxu^\T  &  -(\Bu - \Muw \Mwi \Bt) \end{pmatrix}
$$
and
$$
R = \Mu - \Muw \Mwi \Muwt.
$$

The value function is then equal to:
$$
h_t(\sig{x}(t)) = \sig{x}(t)^\T \sig{n}(t)^\T + \sig{r}(t).
$$

Let the time-varying paraboloid $P = \IVP(\parabol_0)$ diverge in finite time, i.e. $\Tdef{\parabol_0} < \infty$. Since all the touching trajectories are continuous in time, each touching trajectory is defined over $[0,\Tdef{\parabol_0}]$. Their corresponding value function $h$ evaluated along the touching trajectory is as well continuous over $[0,\Tdef{\parabol_0}]$. Therefore, one can prolongates the definition of $P$ until $\Tdef{\parabol_0}$ by continuity of the value function:
$$
P(T)  = \{ z \in \R^{n+1} | \lim_{\substack{t \rightarrow T\\t \leq T}} h(t,z) \leq 0\}.
$$
where $T = \Tdef{\parabol_0}$.
We state this result in the following property:
\begin{prop}
For any $P = \IVP(\parabol_0)$, if the quadratic coefficient of the time-varying paraboloid set $P$ diverges in finite time, then the prolongation to the right of $P$ is defined:
$$
P(T)  = \{ z \in \R^{n+1} | \lim_{\substack{t \rightarrow T\\t \leq T}} h(t,z) \leq 0\}.
$$
\end{prop}




\section{Exact Reachable Set} \label{sec:exact_reachable_set}
\newcommand\defstatefx[3]{
  \expandafter\def\csname Z#1\endcsname{\sig{#2{z}^{#3}}}
  \expandafter\def\csname Z#1s\endcsname{\sig{#2{z}^{#3*}}}
  \expandafter\def\csname x#1\endcsname{\sig{#2{x}^{#3}}}
  \expandafter\def\csname x#1s\endcsname{\sig{#2{x}^{#3*}}}
  \expandafter\def\csname xq#1\endcsname{\sig{#2{x}^{#3}_q}}
  \expandafter\def\csname xq#1s\endcsname{\sig{#2{x}_q^{#3*}}}
  \expandafter\def\csname dxq#1s\endcsname{\dot{\sig{#2{x}}_q^{#3*}}}
  \expandafter\def\csname P#1\endcsname{#2{P}^{#3}}
  \expandafter\def\csname E#1\endcsname{#2{\E}^{#3}}
  \expandafter\def\csname fc#1\endcsname{#2{\fc}^{#3}}
  \expandafter\def\csname fc#1t\endcsname{#2{\fc}^{#3\T}}
  \expandafter\def\csname r#1\endcsname{#2{\r}^{#3}}
  \expandafter\def\csname P#1sup\endcsname{#2{P}^{#3}_+}
  \expandafter\def\csname P#1inf\endcsname{#2{P}^{#3}_-}
  \expandafter\def\csname h#1\endcsname{#2{h}^{#3}}
  \expandafter\def\csname dh#1\endcsname{\dot{#2{h}}^{#3}}
  \expandafter\def\csname w#1s\endcsname{\sig{#2{w}^{#3*}}}
  \expandafter\def\csname G#1\endcsname{\sig{#2{\gamma}^{#3}}}
  \expandafter\def\csname g#1\endcsname{#2{\gamma}^{#3}_0}
  \expandafter\def\csname T#1\endcsname{{#2{T}^{#3}}}
}%
\defstatefx{a}{}{}%
\defstatefx{b}{\tilde}{}%
\defstatefx{c}{}{\prime}%
In this section, we first show that the state constraint \eqref{eq:cons} can be used to redefine the time-varying paraboloids (in Sections \ref{ssec:scaled_paraboloids} and \ref{ssec:domain_def}). Then we define a set of time-varying paraboloid (in Section \ref{ssec:def_parset}). At each time instant, the intersection of these paraboloids is an overapproximation of the reachable set (in Section \ref{ssec:overapprox}). Finally, we prove that when some topological assumption holds about the reachable set, our overapproximation is equal to the reachable set (in Sections \ref{ssec:interior_state}, \ref{ssec:surface_state} and \ref{ssec:exact_reachable_set}).

\subsection{Scaled Paraboloids}\label{ssec:scaled_paraboloids}
Property~\ref{prop:scaled_overapprox} ensures that the overapproximation relationship is still valid when we scale the paraboloid. In this section, a definition of time-varying paraboloid with continuous time scaling is given.

For a given \emph{scaling function} $\Ga \in \Ldloc(\Rp;\Rp)$ (non-negative and locally square integrable function), an initial scaling factor $\ga \geq 1$ and given initial conditions $\parabol_0 = (E_0,f_0,g_0) \in \parabolSet$, as in differential equations (\ref{eq:eq_diff_E},\ref{eq:eq_diff_fc},\ref{eq:eq_diff_g}) in Section~\ref{sec:paraboloid_over_approx}, we can similarly define the initial value problem:
\begin{subequations}
\label{eq:eq_diff_Efcg_gamma}
{%
\fontsize{10}{0} \selectfont
\setlength{\abovedisplayskip}{10pt}
\setlength{\belowdisplayskip}{\abovedisplayskip}
\setlength{\abovedisplayshortskip}{0pt}
\setlength{\belowdisplayshortskip}{3pt}
\begin{align}
\label{eq:diff_E_gamma}
\dE(t) = &
           - \E(t) A - A^\T \E(t) - \Mx 
\\&
           + \pp{\Bt \E(t) + \Mxwt}^\T \Mwi \pp{\Bt \E(t) + \Mxwt}%
\nonumber \\&
           + \sig{\gamma}(t) \E(t)
\nonumber \\
\dfc(t) = &
           - A^\T \fc(t) + (\Mxu + \E(t) \Bu )\sig{u}(t)
\\&
           + (\E(t) B + \Mxw ) \Mwi (\Bt \fc(t)- \Muwt\sig{u}(t))%
\nonumber \\&
+\sig{\gamma}(t) \fc(t)
\nonumber \\
\dr(t) = & \smat{\fc(t)\\\sig{u}(t)}^\T G \smat{\fc(t)\\\sig{u}(t)}%
+\sig{\gamma}(t) \r(t) 
\end{align}}
\end{subequations}
with
\begin{equation}
\label{eq:init_E_gamma}
(\E(0),\fc(0),\r(0)) = (\ga E_0,\ga f_0,\ga g_0).
\end{equation}
The differential equation \eqref{eq:eq_diff_Efcg_gamma} is similar to (\ref{eq:eq_diff_E},\ref{eq:eq_diff_fc},\ref{eq:eq_diff_g}) except that a multiplicator to the constraint is added to the value function and the initial condition of the time-varying paraboloid is scaled by $\ga$.

For $\Ga \in \Ldloc(\Rp;\Rp)$, $\ga \geq 1$, let $P = \IVP(\parabol_0,\ga,\Ga)$ be the time-varying paraboloid with time-varying parameters defined by \eqref{eq:eq_diff_Efcg_gamma} for initial conditions defined by $\parabol_0$.
The Hamiltonian form of the equation can be as well defined and when the quadratic coefficient $\E$ diverges in finite time, the interval of definition of the time-varying paraboloids can be as well prolongated to the closed interval.

The worst disturbance is still expressed by \eqref{def:ws} and Property~\ref{prop:max_dhX} can be restated as:
\begin{prop}
For $\Ga \in \Ldloc(\Rp;\Rp)$, $\ga \geq 1$, and  $\parabol_0 \in \parabolSet$, let $P = \IVP(\parabol_0,\ga,\Ga)$.
For an optimal trajectory $\Zas$ generated by the disturbance $\sig{w}$ defined in \eqref{def:ws} s.t. $\Zas(0) \in \boundary \Pa(0)$, for any $t \geq 0$ it holds:
\begin{equation}
\label{eq:dhtXs_gamma}
\dhtXs(t) = \Ga(t) (\htXs(t)-\xqas(t)).
\end{equation}
\end{prop}
\begin{pf}
Direct derivation from \eqref{eq:eq_diff_Efcg_gamma}.
\end{pf}

When $\Zas(0) \in \boundary \Pa(0)$, the solution to the ODE \eqref{eq:dhtXs_gamma} is, for $t \in \Idef{P}$:
\begin{equation}
\label{eq:sol_htXs}
\htXs(t) =  (1-\ga) \xqas(0) -\lint_0^t \Ga(s) \xqas(s) e^{\int_s^t \Ga(r) dr} ds.
\end{equation}
Since $\int_0^T \gamma(t) \xqas(t) dt$ might not be equal to $0$, trajectories generated by the worst case disturbance $\sig{\ws}$ do not necessarily stay in contact with the time-varying paraboloid and therefore are not touching trajectories. For this reason, we call \textit{optimal trajectories} the trajectories generated by $\sig{\ws}$ given in \eqref{def:ws} .
\begin{prop}
Let $\Zas$ an optimal trajectory of $\Pa$ s.t. $\Zas(0) \in \boundary \Pa(0)$, if:
$$(1-\ga) \xqas(0) - \lint_0^t \Ga(\tau) \xqas(\tau) d\tau = 0$$ for any $t \geq 0$ and $\htXs(0) = 0$, then $\Zas$ is a touching trajectory of $\Pa$.
\end{prop}
\begin{pf}
Using \eqref{eq:sol_htXs}.
\end{pf}

Therefore, for any other trajectory of the constrained system $\sys$, $h$ is a decreasing function of time along the trajectory.
In this case, Theorem~\ref{thm:overapproximation} can be rewritten for continuously scaled time-varying paraboloid:
\begin{thm}
\label{thm:overapprox_scaled}
For a set of initial states $\parabol_0$, a scaling function $\Ga \in \Ldloc(\Rp; \Rp)$ and an initial scaling factor $\ga \geq 1$, let $P = \IVP(\parabol_0,\ga, \Ga)$. The reachable set $\reach(\parabol_0,t)$ of $\sys(\parabol_0,t)$, $t>0$, is overapproximated by $P(t)$, i.e.:
$$
\forall t \in \Idef{P}, \reach(\parabol_0,t) \subseteq P(t) \cap \Zp.
$$
\end{thm}
\begin{pf}
By integration of \eqref{eq:sol_htXs} over the interval $[0,t]$.
\end{pf}

\subsection{Definition domain of time-varying paraboloid}\label{ssec:domain_def}
\label{sec:unstablecase}
\newcommand{\IVPk}{\IVP_\kappa}
In the case where the DRE in \eqref{eq:eq_diff_E} does not have any convergent solution for any positive definite initial condition (i.e. when system defined in Section~\ref{sec:problem_statement} is unstable), works based on the LQR formulation of the IQC cannot overapproximate the reachable set for any $t \geq 0$ (for any initial condition of the DRE, the solution to the DRE is defined over a finite escape time). In this part, we show that for any system (whether it is stable or unstable) there is always a $\Ga \in \Ldloc(\Rp;\Rp)$ such that the corresponding time-varying paraboloid is defined over $\Rp$. 
This is one of the key advantages of our approach.

\newcommand{\pE}{E_0}
For a given positive definite initial condition $E_0 \mgeq 0$, for the scaling factor $$\Ga(\cdot) = \kappa \geq 0$$ over $\Rp$ and a scaling factor $\ga = 1$,  the corresponding solution $\E$ to the DRE \eqref{eq:diff_E_gamma} does not diverge over $\Rp$ if:
\begin{equation}
\label{eq:stab_care_k}
\begin{bmatrix}
-\pE A - A^\T \pE -\Mx + \kappa \pE& \Bt \pE + \Mxwt\\
(\Bt \pE + \Mxwt)^\T & -\Mw
\end{bmatrix}
\mge 0.
\end{equation}
Since $-\Mw \mge 0$, the Schur complement of \eqref{eq:stab_care_k} leads to the equivalent non-negativity condition:
$$
\overline{E}_0 + \kappa \pE \mgeq 0 
$$
where 
$$
\begin{aligned}
\overline{E}_0 = 
&- \pE A - A^\T \pE - \Mx +\\
& \pp{\Bt \pE + \Mxwt}^\T \Mwi \pp{\Bt \pE + \Mxwt}.
\end{aligned}
$$
By choosing $\kappa$ such that
$$
\kappa > \frac{\norm{\overline{E}_0}}{\norm{\pE}} ,
$$
then there is a convergent solution $\E$ to the DRE \eqref{eq:diff_E_gamma}. Therefore, for any given $\pE \mge 0$, there exists a $\kappa>0$ such that \eqref{eq:stab_care_k} is satisfied.

\begin{prop}
\label{prop:cvg_tvp}
There is a $\Ga \in \Ldloc(\Rp;\Rp)$ such that $P = \IVP(\parabol_0,1,\Ga)$ is defined over $\Rp$.
\end{prop}

By Property~\ref{prop:scaled_overapprox}, $\Pa(t) \cap \Zp$ is an overapproximation of the reachable set $\reach(\parabol_0,t)$ for any $t \geq 0$.

\begin{rem}
At a given time, the reachable set of an IQC system is always bounded.
Whether the DRE \eqref{eq:eq_diff_E} has a solution or not over a given interval $[0,T]$, $T>0$, we can always bound the set of reachable states. This result is the main difference with other works in reachable set overapproximation for IQC systems (see \cite{savkin1995recursive,guseinov2009approximation,gusev2017extremal,jonsson2002robustness,seiler2019finite}).
Up to the knowledge of the authors, other works reachable set overapproximation for IQC systems only use results from Theorem~\ref{thm:overapproximation}. These results are dependent on the existence of a solution to the DRE. By taking into account the constraint over the interval of integration and not only at the terminal time, we have less conservative results.
\end{rem}

\begin{figure*}[t]
\newcommand{\imgLarge}{0.22\textwidth}
\includegraphics[angle=90,height=\imgLarge]{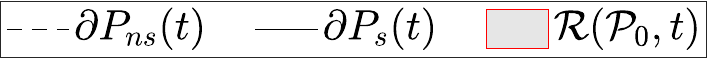}
\includegraphics[width=\imgLarge]{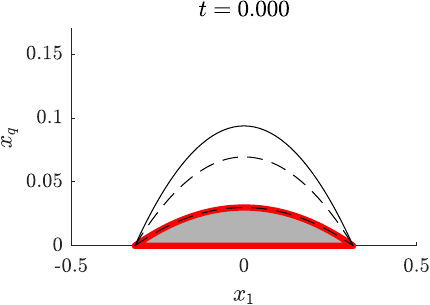}
\includegraphics[width=\imgLarge]{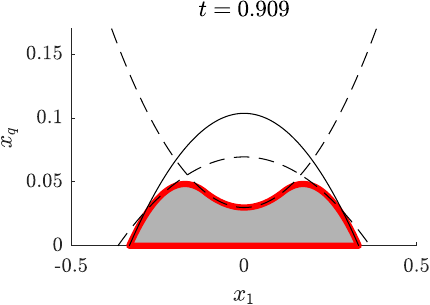}
\includegraphics[width=\imgLarge]{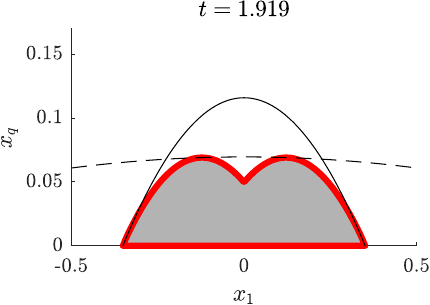}
\includegraphics[width=\imgLarge]{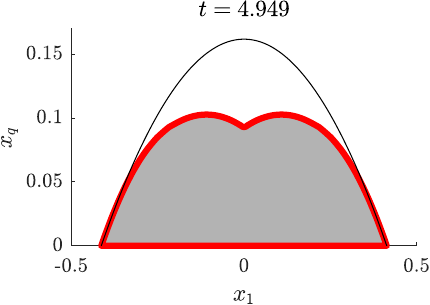}
\caption{\label{fig:unstable_system}Stabilized constraint $P_s(t)$ versus finite escape time constraints $P_{ns}(t)$. The time-varying paraboloid $P_{ns}$ is defined over $[0,3]$ whereas $P_s$ is defined over $\Rp$.}
\end{figure*}

\begin{exmp}
Figure~\ref{fig:unstable_system} shows plots of the reachable set of the unstable system $\sys(\parabol_0)$ defined by parameters:
\[
A = -1, \,
B = 1, \,
B_u = 0, \,
M = \smat{1&0&0\\0&1&0\\0&0&-0.9}
\]
and a zero input signal $u$.
The set of initial states $\Z_0$ is a paraboloid $\Z_0 = \parabol_0 = \parabol(\E_0,\fc_0,\r_0)$ with:
\[
\E_0 = 1, \,
\fc_0 = 0 \textrm{ and }
\r_0 = 0.015.
\]
The solution to DRE~\eqref{eq:diff_E_gamma} for $\Ga = 0$ and $\ga = 0$ has a finite escape time and diverges at $\Tdef{P_{ns}} = 1.7$. 
The solution to DRE~\eqref{eq:diff_E_gamma} for $\Ga = 0$ and $\ga = 1$ is defined over $\Rp$. 
\end{exmp}

\subsection{Overapproximation with an intersection of time-varying paraboloids}\label{ssec:def_parset}
In this section, a set of time-varying paraboloids is defined. At a given time, the intersection of the paraboloids gives better overapproximations of the reachable set. With additional assumptions about the topology of the reachable set, the reachable set is exactly characterized. This approach relies on the use of Property~\ref{thm:overapproximation} and preliminary results showing that for any state of the overapproximation, there exists a trajectory in $\sys(\parabol_0,t)$, $t>0$, leading to this state.

Let $\parset$ be defined as follows:
\begin{equation}
\label{def:exact_reach_set}
\parset = \{ \IVP(\parabol_0,\ga,\Ga) | \Ga \in \Ldloc(\Rp;\Rp), \Ga \geq 0, \ga \in \R, \ga \geq 1\}.
\end{equation}
$\parset$ corresponds to the set of all time-varying paraboloids with initial conditions $\parabol_0$ and generated by the set of non-negative scalings $\Ga \in \Ldloc(\Rp;\Rp)$ and the set of initial scaling factors $\ga \geq 1$.
Let 
\begin{equation}
\label{def:parset}
\parset(t) = \left\{ P \in \parset \middle| t \in \Idef{P} \right\}
\end{equation}
the set of all the defined time-varying paraboloids at time $t \geq 0$.
By Property~\ref{prop:cvg_tvp}, $\parset(t) \neq \emptyset$ for any $t \geq 0$. In other words, $\parset$ is defined over $\Rp$ as well.

\begin{figure}[ht!]
\centering
\includegraphics[width=\columnwidth]{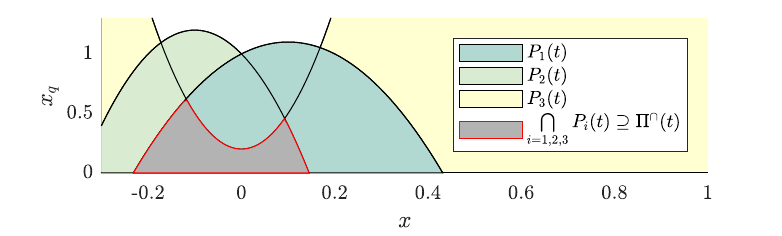}
\caption{\label{fig:describe_Pi}For a given $t \geq 0$, let 3 time-varying paraboloids $P_i \in \parset$, $i=1,2,3$. Light color shaded area are their corresponding parabolic set $P_i(t)$ at $t$, $i=1,2,3$. Grey color shaded is their intersection. By \eqref{def:parinter}, $\parinter(t)$ is a subset of $P_1(t) \cap P_2(t) \cap P_3(t)$.}
\end{figure}

For $t \geq 0$, let 
\begin{equation}
\label{def:parinter}
\parinter(t) = \bigcap_{P \in \parset(t)} P(t)
\end{equation}
the intersection of all the defined time-varying paraboloids $P$ of $\parset$ at time $t$ (see Figure~\ref{fig:describe_Pi}). Since $\parinter(\cdot)$ is defined over $\Rp$, $\parinter(\cdot)$ is defined over $\Rp$.

\bigskip 
We now prove that, when some assumptions about the topology of $\parinter$ hold (Assumption~\ref{hyp:boundary_reached} and \ref{hyp:dvg_paraboloid_regular}), we have $\reach(\parabol_0,t) = \parinter(t) \cap \Zp$, for any $t \geq 0$ (Theorem~\ref{thm:exact_reachable_set}, Section~\ref{ssec:exact_reachable_set}). To achieve that:
\begin{itemize}
\item we prove the overapproximation relationship $\reach(\parabol_0,t) \subseteq \parinter(t)$ (Section~\ref{ssec:overapprox});
\item we prove that any state $\extstate{x}{x_q} \in \parinter(t)$ is reachable from a state $\extstate{x}{x_q'} \in \boundary \parinter(t)$ with $x_q \leq x_q'$ (Section~\ref{ssec:interior_state});
\item for a state $z_t \in \boundary \parinter(t)$, we find a touching trajectory $\sig{z^*} = \extstate{\sig{x^*}}{\sig{x_q^*}}$ of $\parinter$ such that $\sig{z^*}(t) = z_t$. This touching trajectory $\extstate{\sig{x^*}}{\sig{x_q^*}}$ of $\parinter$ satisfies the state constraint $\sig{x_q}(\cdot) \geq 0$ over $[0,t]$  (Section~\ref{ssec:surface_state});
\item finally, we conclude that any $z_t \in \parinter(t)$ is reachable from $\parabol_0$, thus $\reach(\parabol_0,t) = \parinter(t) \cap \Zp$ (Section~\ref{ssec:exact_reachable_set}).
\end{itemize}


\subsection{Overapproximation Relationship} \label{ssec:overapprox}
Theorem~\ref{thm:overapprox_scaled} states that each time-varying paraboloid defined in Section~\ref{ssec:scaled_paraboloids} is an overapproximation of the reachable set. An intersection of many time-varying paraboloids is as well an overapproximation of the reachable set. 
\begin{prop}
\label{prop:overapprox_inter}
$\reach(\parabol_0,t) \subseteq \parinter(t) \cap \Zp$ for any $t \geq 0$.
\end{prop}
\begin{pf}
This is a direct consequence of Theorem~\ref{thm:overapprox_scaled} and \eqref{def:parinter}.
\end{pf}

\begin{figure*}[t]
\newcommand{\imgLarge}{0.23\textwidth}
\centering
\includegraphics[angle=90,scale=0.4]{img/dummy_legend.pdf}
\includegraphics[width=\imgLarge]{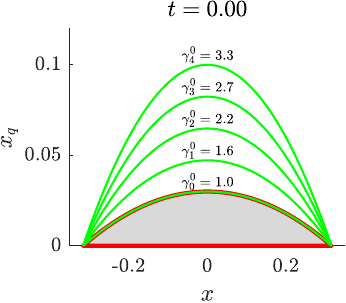}
\includegraphics[width=\imgLarge]{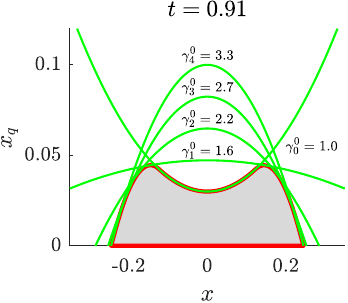}
\includegraphics[width=\imgLarge]{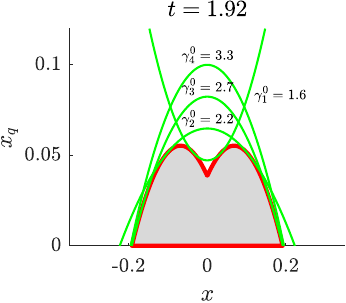}
\includegraphics[width=\imgLarge]{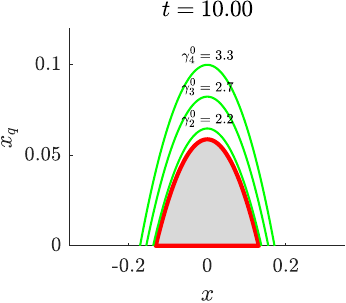}
\caption{\label{fig:energy_levels}%
Time-varying paraboloids overapproximating the reachable set at different time instants $t$ in $\{0.00,0.91,1.62,10.00\}$ for different scalings. Scaling functions (i.e. time-varying scaling factors, see Section~\ref{ssec:scaled_paraboloids}) are equal to zeros, $\Ga_i = \sig{0}$, and initial scaling factors $\scaling_i^0 \geq 1$ are respectively equal to $1.0$, $1.6$, $2.2$, $2.7$ and $3.3$ for $i=0,\dots,4$. The shaded regions are the reachable set $\reach(\parabol_0,t)$, the thin lines are the boundary of the overapproximation $P(t)$ of Theorem~\ref{thm:overapproximation}.}
\end{figure*}

\begin{exmp}[Continued from Example~\ref{ex:1dsys}] \label{ex:1dsys_cont}
In the case where the solution to \eqref{eq:eq_diff_E} does not converge (i.e. $\Ei_0<\E^-$), Figure~\ref{fig:energy_levels} shows several paraboloid trajectories with different initial scaling factors. Scaling functions are equal to 0 and initial scaling factors $\scaling_i$ are greater than $1$, $\parabol_0 \cap \Zp \subset \scaling_i \parabol_0 \cap \Zp$. Therefore, each time-varying paraboloid is a valid constraint that bounds $\reach(\parabol_0,t)$, $t \in \Idef{\parset}$ (Theorem~\ref{thm:overapproximation}). Therefore, $\reach(\parabol_0,t) \subseteq P^\cap(t) = P_0(t) \cap P_1(t) \cap \dots \cap P_4(t)$ where $P_i = \IVP(\parabol_0,\scaling_i , \sig{0})$,  and $\scaling_i$ are resp. equal to $1$, $1.6$, $2.2$, $2.7$ and $3.3$ for $i = 0,\dots,4$. In this case, the overapproximation $P^\cap(t)$ is strictly included in $P_{0}(t)$.
\end{exmp}
Observations in Example~\ref{ex:1dsys_cont} motivate the use of multiple time-varying paraboloids to get better overapproximations of the reachable set $\reach(\parabol_0,t)$, $t>0$.

\subsection{Past trajectory for states in the overapproximation's interior} \label{ssec:interior_state}
Property~\ref{prop:consume_energy} is already presented in \cite[][Property~7]{rousse:hal-02049158}, we restate it hereby for completeness.

\newcommand{\ti}{{t^*}}%
\newcommand{\tf}{t}%
\newcommand{\Ze}{\Z_\ee}%
\newcommand{\g}{r}%
\newcommand{\dg}{\dot{r}}%
\newcommand{\xs}{\sig{x^*}}%
\newcommand{\xqs}{\sig{\xq^*}}%
\newcommand{\DE}{\Delta}
\newcommand{\Dfc}{\mu}
\newcommand{\ee}{\epsilon}
\newcommand{\dd}{\delta}

Property~\ref{prop:consume_energy} shows that the state $\extstate{x}{\alpha x_q}$ is reachable from the given state $\extstate{x}{x_q}$ for any given $\alpha \in [0,1]$.
\begin{prop}
\label{prop:consume_energy}
For $t \geq 0$, if $\extstate{x}{x_q} \in \reach(\parabol_0,t)$ then $\extstate{x}{\alpha x_q} \in \reach(\parabol_0, t)$ for all $\alpha \in [0,1]$.
\end{prop}
\begin{pf}
Let $f: t,x \mapsto A x + B \sig{w}(t) + B_u \sig{u}(t)$. Since for any $(t,x) \in \Rp \times \R^n$, $f(.,x)$ is locally measurable over $\Rp$, $f(t,.)$ is Lipschitz over $\R^n$, \eqref{eq:ode} has a unique solution $\sig{x}$ (see \cite{schuricht2000ordinary}, Theorem~1.1) that is time-continuous.
Therefore, for a trajectory $\extstate{\sig{x}}{\sig{x_q}} \in \sys(\parabol_0,T)$, $T>0$, $\sig{x}$ is time-continuous.%

For $\ee>0$, let $\sig{w} \in \Ld([0,t+\ee];\R^m)$, s.t. $\sig{w}^\T(s) \Mw \sig{w}(s) = -(1-\alpha) \sig{x_q}(t) \frac{1}{\ee}$ when $s \in [t,t+\ee]$. Then $$\int_t^{t+\ee} \sig{w}^\T(s) \Mw \sig{w}(s) ds \rightarrow -(1-\alpha) \sig{x_q}(t)$$ when $\ee \rightarrow 0$. Using Cauchy-Schwartz inequality:
$$
\abs{\lint_t^{t+\ee} (-\Mw)^{\frac{1}{2}} \sig{w}(s) ds}
\leq
\sqrt{\ee}
\sqrt{\lint_t^{t+\ee} -\sig{w}^T(s) \Mw \sig{w}(s) ds}
$$
and the time-continuity of $\sig{x}$, the quantity
$$
\lint_t^{t+\ee} \smat{\sig{x}(s)\\\sig{u}(s)\\0}^\T M \smat{\sig{x}(s)\\\sig{u}(s)\\\sig{w}(s)} ds \rightarrow 0$$
when $\ee \rightarrow 0$. By integration, $\sig{x}_q(t+\ee) \rightarrow \alpha \sig{x_q}(t)$ when $\ee \rightarrow 0$. Since $\sig{x}$ is time-continuous, $\sig{x}(t+\ee) \rightarrow \sig{x}(t)$ when $\ee \rightarrow 0$.
By continuity of $\sig{u}$, $\sig{x}$ and $\sig{w}$ over $[t,t+\ee]$, $\sig{x_q}$ is continuous over $[t,t+\ee]$. Then, there exists a $t' \in [t,t+\ee]$ such that $\sig{x_q}(\tau) \geq \alpha \sig{x_q}(t) \geq 0$ for all $\tau \in [t,t']$ and $\sig{x_q}(t') \rightarrow \alpha \sig{x_q}(t)$ when $\ee \rightarrow 0$. Therefore, the constraint $\sig{x_q}(\cdot)\geq 0$ is satisfied over $[t,t']$ and the trajectory $(\sig{x},\sig{x_q})$ is a valid trajectory of $\sys(\parabol_0,t')$ for all $t \leq t'$.
\end{pf}

\subsection{Past trajectory for states in the overapproximation's boundary} \label{ssec:surface_state}%
\renewcommand{\ll}{\lambda}%
\newcommand{\dde}{\sig{\delta}}%
\newcommand{\ddez}{\delta_0}%
In this section, touching trajectories of $\parinter$ are identified. We show that all these touching trajectories satisfies the state constraint \eqref{eq:state_cons}.

The value function $\hb$ of a time-varying paraboloid $\Pb \in \parset$ can be approximated at the first order along a touching trajectory $\sig{z^*}$ of another time-varying paraboloid $\Pa \in \parset$ when their scaling functions $\Gb$ and $\Ga$ and initial scaling factor $\gb$ and $\ga$ are close. In this part, we compute this first order approximation when $\Gb = \Ga + \dde$ and $\gb = \ga + \ddez$ for small variations $\dde \in \Ld(\Rp;\R)$ and $\ddez \in \R$  (i.e. when $\norm{\dde} + \abs{\ddez}$ tends to $0$).

To prove that the reachable set $\reach(\parabol_0,t)$, $t>0$, is exactly described by $\parinter(t)$, we show that for any optimal trajectory $\Zas$ of $\Pa \in \parinter$ s.t. $\Zas$ is violating the constraint \eqref{eq:state_cons}, there is a $\Pb \in \parset$ such that the end point $\Zas(t)$ does not belong to $\Pb(t)$ and therefore to $\parinter(t)$. To do so, we will study the value function of a time-varying paraboloid $\Pb$ for touching trajectories of $\Pa$. 

\newcommand{\la}{\lambda}
\newcommand{\lb}{\tilde{\lambda}}
\begin{prop}
\label{prop:contact_traj_first_order}
For $\Ga \in \Ldloc(\Rp,\Rp)$ and $\ga \geq 1$, let the corresponding time varying paraboloid $\Pa = \IVP(\parabol_0,\ga, \Ga)$.
For any $t$ in the open set of $\Idef{\Pa}$, it exists $\epsilon>0$ and $H>0$ s.t. for any $\dde \in \Ld(\Rp;\R)$, $\norm{\dde} \leq \epsilon$,
for any $\ddez \in \R$, $\abs{\ddez} \leq \epsilon$, 
 $\Pb = \IVP(\parabol_0,\gb, \Gb)$ where $\Gb  = \Ga + \dde$ and $\gb = \ga + \ddez$ s.t. $t$ belongs to the open set of $\Idef{\Pb}$, 
let $\hb_t$ the value function of $\Pb(t)$ and $\Zas = \extstate{\xas}{\xqas}$  an optimal trajectory of $\Pa$, it holds:
\begin{equation*}
\abs{ \hb_t(\Zas(t)) - \beta(t) }  \leq H  \epsilon^2
\end{equation*}
where
\begin{equation*}
\beta(t) = \frac{\ddez}{\ga} \xqas(0)+\lint_0^t \dde(s) \psi(s) e^{\int_s^t (\Ga(r) + \dde(r))dr}ds
\end{equation*}
and
\begin{equation*}
\psi(s) = \xqas(s)-\lint_s^t \Ga(\tau) \xqas(\tau) e^{\int_\tau^s \Ga(r) dr} d\tau.
\end{equation*}
\end{prop}
\begin{pf}
Let $(\Ea,\fca,\ra)$ and $(\Eb,\fcb,\rb)$ (resp.) be parameters of $\Pa$ and $\Pb$(resp.), and
$$
\sig{\nu} = (\Ea-\Eb) \sig{\xas} - (\fca - \fcb).
$$
Using (\ref{eq:ode}, \ref{def:ws}, \ref{eq:eq_diff_Efcg_gamma}), $\sig{\nu}$ satisfies the linear time varying differential equation:
\begin{equation}
\label{eq:nu_ode}
\dot{\sig{\nu}}(\tau) = A_\nu(\tau) \sig{\nu}(\tau) - \dde(\tau) \sig{n}(\tau).
\end{equation}
with $\sig{n} = \Ea \xas -\fc$ and $A_\nu(\tau) = -A^\T + \Mxw \Mwi \Bt + \Eb(\tau) B \Mwi \Bt + \Ga(\tau) I$. By \eqref{eq:init_E_gamma}, initial values of $\Pa$ and $\Pb$ satisfies 
$$
\frac{1}{\ga} \Pa(0) = \frac{1}{\gb}\Pb(0) = \parabol_0,$$
therefore $\sig{\nu}(0) = \ddez (E_0 \xas(0) - f_0)$ where $(E_0,f_0,g_0) = \parabol_0$.
Since $t$ belongs to the open set of $\Idef{\Pb}$, $\Eb(\cdot)$ is bounded over $[0,t]$ (the discontinuity of $\Eb$ can only occur at the final integration time). By time-continuity of $\Eb(\cdot)$ over $[0,t]$, there is a scalar $K>0$ that bounds $\norm{\Eb(\cdot)}$ over $[0,t]$. Then, since $\Ga$ is measurable, there exists a measurable function $L \in \Ld([0,t];\Rp)$ such that:
\begin{equation}
\label{eq:nu_diff_eq}
\norm{A_\nu(\tau)} \leq L(\tau)
\end{equation}
over $\tau \in [0,t]$.
We can integrate \eqref{eq:nu_diff_eq}:
$$
\norm{\sig{\nu}(t)-\sig{\nu}(0)} \leq \lint_0^t L(\tau) \norm{\sig{\nu}(\tau)} d\tau + \lint_0^t \abs{\dde(\tau)} \norm{\sig{n}(\tau)} d\tau.
$$ 
Since $l \mapsto \int_0^l \abs{\dde(\tau)} \norm{\sig{n}(\tau)} d\tau$ is a non decreasing function over $[0,t]$, by applying the Gr\"onwall inequality, we get:
\begin{equation}
\label{eq:bound_nu}
\norm{\sig{\nu}(t)} \leq \left( \norm{\sig{\nu}(0)} + \lint_0^t \abs{\dde(\tau)} \norm{\sig{n}(\tau)} d\tau \right)
\,
e^{\int_0^t L(\tau) d\tau} .
\end{equation}%
\newcommand{\diffh}{q}%
Let
\begin{equation}
\label{eq:def_q}
\diffh(t)  = \ha_t(\Zas(t)) - \hb_t(\Zas(t)).
\end{equation}
$\Zas$ is an optimal trajectory of $\Pa$ s.t. $\Zas(t) \in \boundary \Pa(t)$, therefore, $\ha_t(\Zas(t)) = 0$, therefore, using \eqref{eq:dhtXs_gamma}:
\begin{equation}
\label{eq:h_eq_int}
\ha_\tau(\Zas(\tau)) = \lint_\tau^t \Ga(s) \xqas(s) e^{\int_\tau^s \Ga(r) dr} ds.
\end{equation}
Using (\ref{eq:ode}, \ref{def:ws}, \ref{eq:eq_diff_Efcg_gamma}, \ref{eq:h_eq_int}), $\diffh$ satisfies:
$$
\begin{aligned}
\dot{\diffh}(\tau) =
& -\sig{\nu}^\T(\tau)  B \Mwi \Bt \sig{\nu}(\tau) + \Ga(\tau) \ha_\tau(\Zas(\tau)) \\
&- \Gb(\tau) \hb_\tau(\Zas(\tau)) - \Ga(\tau) \xqas(\tau) + \Gb(\tau) \xqas(\tau) .
\end{aligned}
$$
Using \eqref{eq:def_q} and \eqref{eq:h_eq_int}: 
\begin{equation}
\label{eq:ode_q}
\begin{aligned}
\dot{\diffh}(\tau) =
& -\sig{\nu}^\T(\tau)  B \Mwi \Bt \sig{\nu}(\tau) \\
&+ \dde(\tau) \psi(\tau) + (\Ga(\tau)+\dde(\tau)) \diffh(\tau)
\end{aligned}
\end{equation}
where 
$$
\psi(\tau) = \xqas(\tau) - \lint_\tau^t \Ga(s) \xqas(s) e^{\int_\tau^s \Ga(r) dr} ds
$$
with initial condition $\diffh(0) = \ha_0(\Zas(0))-\hb_0(\Zas(0))$. Since $\Zas$ is a touching trajectory of $\Pa$, it holds $\ha_0(\Zas(0)) = 0$, therefore:
$$
\xqas(0) = -\ga ( \xas(0)^\T E_0 \xas(0) - 2 f_0^\T \xas(0) + g_0).
$$
Therefore, $\hb_0(\Zas(0))$ satisfies:
$$
\hb_0(\Zas(0)) = -\frac{\ddez}{\ga} \xqas(0)
$$
and $\diffh(0) = \frac{\ddez}{\ga} \xqas(0)$.

Since $t$ belongs to the open set of $\Idef{\Pb}$, the optimal trajectory $\Zas$ and $\sig{\nu}$ are defined and continuous over $[0,t]$. Moreover, since $\Ga$ and $\dde$ are measurable over $[0,t]$, the solution to the linear time-varying equation \eqref{eq:ode_q} exists over $[0,t]$ and is:
$$
\begin{aligned}
\diffh(\tau) =& -\hb_0(\Zas(0))+
\lint_0^\tau \bigg[ \Big( -\sig{\nu}(s)^\T B \Mwi \Bt \sig{\nu}(s)
\\& + \dde(s) \psi(s)\Big) e^{\int_s^t (\Ga(r)+\dde(r)) dr} \bigg]ds.
\end{aligned}
$$
Then, using \eqref{eq:bound_nu}:
$$
\abs{
\diffh(t) - \frac{\ddez}{\ga} \xqas(0) -
\lint_0^t
\dde(s) \psi(s) e^{\int_s^t (\Ga(r)+\dde(r)) dr} 
ds
}
\leq
H  \epsilon^2
$$
with
\begin{equation}
\label{eq:def_H}
H = t R K^2 N (\norm{\sig{n}(0)}^2 + \norm{\sig{n}}^2)
\end{equation}
a finite constant
where $R = \norm{B \Mwi \Bt}$, $K = \exp{\int_0^t L(\tau) d\tau}$ and
$N = \int_0^t e^{\int_s^t (\Ga(r)+\dde(r)) dr} ds$.
This ends the proof.
\end{pf}

\begin{rem}
When $\Ga = \sig{0}$, Property~\ref{prop:contact_traj_first_order} matches with \cite[Property~9]{rousse:hal-02049158} where their initial scaling correspond to our initial scaling factor $\ga$.
\end{rem}

Property~\ref{prop:reject} gives conditions where the sign of $\hb_\tf(\Zas(t))$ is only determined by its first order approximation defined in Property~\ref{prop:contact_traj_first_order}.

\begin{prop}
\label{prop:reject}
Let $\sig{z^*}$ a touching trajectory of $\Pa = \IVP(\parabol_0,\ga,\Ga)$ for $\Ga \in \Ldloc(\Rp;\Rp)$, $\ga \geq 1$ given and $t \in \Idef{\Pa}$ given. If there is a $\dde \in \Ldloc(\Rp;\R)$ and a $\ddez \in \R$, s.t. $\norm{\dde} \leq \epsilon$ and $\abs{\ddez} \leq \epsilon$
and $t \in \Idef{\Pb}$ (where $\Pb = \IVP(\parabol_0,\ga+\ddez, \Ga+ \dde)$) and
\begin{equation}
\label{eq:first_order_cond}
H \epsilon^2 \le \abs{
\frac{\ddez}{\ga} \xqas(0)
+\lint_0^t
\left[
	\dde(s) 
	\psi(s)
	e^{\int_s^t (\Ga(r)+\dde(r)) dr} 
\right]ds
},
\end{equation}
then the sign of
\begin{equation*}
-\frac{\ddez}{\ga} \xqas(0)
-\lint_0^t
\left[
	\dde(s) 
	\psi(s)
	e^{\int_s^t (\Ga(r)+\dde(r)) dr} 
\right]ds
\end{equation*}
is equal to the sign of $\hb_\tf(\Zas(t))$
where $\hb_t$ is the value function of $\Pb(t)$ and $H>0$ defined in \eqref{eq:def_H} and
$$
\psi(s)=
    \xqas(s)
    -\lint_s^t \Ga(\tau) \xqas(\tau) e^{\int_\tau^s \Ga(r) dr} d\tau 
.
$$
\end{prop}
\begin{pf}
This is a direct consequence of Property~\ref{prop:contact_traj_first_order} and of the property: $(\abs{a-b} \leq c) \land (c < \abs{b}) \Rightarrow \sign(a)=\sign(b)$ for $a,b,c \in \R$.
\end{pf}

Provided the existence of a $(\dde,\ddez) \in \Ld(\Rp;\R) \times \R$ such that $\Ga+\dde \geq 0$ and $\ga + \ddez \geq 1$, the first order approximation of the value function of $\Pb = \IVP(\parabol_0, \ga +\ddez, \Ga+\dde)$ gives a way to identify time varying paraboloids $\Pb$ that belongs to $\parset$ such that an invalid trajectory with an end state $z_t \in \boundary P(t)$ (meaning with initial state outside of the initial set $\parabol_0$ or a trajectory violating the constraint) does not belongs to $\Pb(t)$ and therefore, does not belong to $\parinter(t)$.

Property~\ref{prop:case_psi} states that the touching trajectories of $\parinter$ satisfy the state constraint \eqref{eq:state_cons}.
Property~\ref{prop:case_psi} is proven by choosing a valid trajectory candidate. If this trajectory violates the state constraint \eqref{eq:state_cons}, then  Property~\ref{prop:reject} provides a proof that this trajectory does not belongs to the overapproximation $\parinter$.
\begin{prop}
\label{prop:case_psi}
For $\Pa \in \parinter$, if $z_t \in \boundary \parinter(t)$ and $z_t \in \boundary \Pa(t)$ for $t$ in the open set of $\Idef{\Pa}$, then the optimal trajectory $\Zas$ of $\Pa$ such that $\Zas(t) = z_t$ is a valid touching trajectory of $P$ and $\parinter$.
\end{prop}
\begin{pf}
Let $\psi: \Rp \mapsto \R$ defined for $s \geq 0$ by:
$$
\psi(s) = \xqas(s)-\lint_s^t \Ga(\tau) \xqas(\tau) e^{\int_\tau^s \Ga(r) dr} d\tau
.$$
Let $\tau \in [0, t]$ and $I = [\tau, t]$.
\begin{itemize}[topsep=0pt,itemsep=-1ex,partopsep=1ex,parsep=1ex]
\item \emph{Case 1, $\xqas(0)<0$}: with $\ddez>0$, using Property~\ref{prop:reject}, $\Zas(t) \notin \Pb(t)$ where $\Pb \in \parset$ since $\ddez+ \ga \geq 1$, so $\Zas(t) \notin \parinter(t)$.
\item \emph{Case 2, $\psi(\cdot) < 0$ over $I$}: any $\dde(\cdot) \geq 0$ over $I$ and $\dde(\cdot) = 0$ elsewhere such that $\int_0^t \dde(s)\psi(s) ds \neq 0$ and for $\ddez = 0$, using Property~\ref{prop:reject}, $\Zas(t) \notin \Pb(t)$ where $ \Pb \in \parset$ since $\Ga + \dde \geq 0$, so $\Zas(t) \notin \parinter(t)$.
\item \emph{Case 3, $\psi(\cdot) > 0$ over the open of $I$ and there is a $l\in I$, s.t. $\int_{s\in [\tau,l]} \Ga(s) \xqas(s)ds \neq 0$}: since $\Ga \geq 0$, there exists a $\dde \leq 0$ such that $\Ga + \dde \geq 0$ and for $\ddez = 0$, using Property~\ref{prop:reject}, $\Zas(t) \notin \parinter(t)$.
\item \emph{Case 4, $\psi(\cdot) = 0$ over $I$}:
since $\xqas$ is continuous over $\Rp$, and since $\Ga$ is locally square integrable, $\psi(t) = 0 \Rightarrow \xqas(t) = 0$, therefore $\xqas(\cdot) = 0$ over $I$. Consequently $\ha_\tau(\Zas(\tau)) = 0$ for $\tau \in I$.
\end{itemize}
Cases 1 to 4 show that for $\Zas(t) \in \boundary \parinter(t)$:
\begin{itemize}[topsep=0pt,itemsep=-1ex,partopsep=1ex,parsep=1ex]
\item either $\forall l \in I, \int_\tau^l \Ga(\tau) \xqas(\tau) d\tau=0$ and $\psi(\cdot)= \xqas(s) >0$;
\item nor $\xqas(l)=0$ for $l \in I$. 
\end{itemize}
Let a partition $[0,t] = \bigcup_{i\in \N} I_i$ such that over each open interval $I_i$, $\psi(\cdot) \bowtie_i 0$ with $\bowtie_i \in \{<,>,=\}$.
We deduce that for any $s \in [0,t]$:
$\xqas(s) \geq 0$ and $\int_I \Ga(\tau) \xqas(\tau)d\tau = 0$.
$\Zas$ is a valid trajectory, i.e. the constraint \eqref{eq:state_cons} is satisfied.
By \ref{eq:dhtXs_gamma}, $\Zas$ is a touching trajectory of $P$. Moreover, since $\Zas(0) \in \parabol_0$, $\Zas$ is as well a touching trajectory of $\parinter$.
%
%
%
\end{pf}

\bigskip
Since $\parinter(t)$ is an intersection of closed sets, $\parinter(t)$ is closed as well. In the general case, for an infinite intersection $\mathcal{Y}^\cap = \bigcap_{i\in \N} Y_i$ of closed sets $Y_i$, $i \in \N$, any boundary point $y \in \boundary \mathcal{Y}^\cap$ does not necessarily belongs to the boundary of any $Y_i$, $i\in \N$ (e.g. $\bigcap_{\ee \in ]1,2]} [-\ee,\ee] = [-1,1]$, but there is no $\ee \in ]1,2]$ such that $1 \in \boundary [-\ee,\ee]$).
The following assumption states that for every state on the boundary of the overapproximation $\parinter(t)$, $t>0$, there exists a time-varying paraboloid $\Pa$ such that this state belongs as well to the boundary of the $\Pa(t)$.
\begin{assum}
\label{hyp:boundary_reached}
For any $z_t \in \boundary \parinter(t)$, there is a $\Pa \in \parset$ such that $z_t \in \boundary P(t)$.
\end{assum}
This assumption is not a strong one and has been proved for simpler cases (see \cite[][Property~11]{rousse:hal-02049158}).

\bigskip
In Property~\ref{prop:case_psi}, the existence of $\Gb$ and $\ga$ is conditioned by $t$ belonging to the open domain $\Idef{\Pb}$; to ensure this, $\norm{\E(\cdot)}$ is assumed to be bounded over $[0,T]$ (by considering the case where $t$ is in the open set of $\Idef{\Pa}$).
In the general case, the boundedness of $\norm{\E(\cdot)}$ is not granted (see the unstable case in Example~\ref{ex:1dsys} and Figure~\ref{fig:bounding_paraboloid_low_energy}).
Assumption~\ref{hyp:dvg_paraboloid_regular} states that for any state on the boundary of the overapproximation $\parinter(t)$, $t>0$, there is neighbor state on the boundary of $\Pb(t)$ where $\Pb$ is a time-varying paraboloid of $\parset$ not diverging at $t$ (i.e. $t$ belongs to the interior of $\Tdef{P}$).

\begin{assum}
\label{hyp:dvg_paraboloid_regular}
For $t>0$, for all $\epsilon > 0$, for any $z_t \in \boundary \parinter(t)$ such that $z_t \in \boundary \Pa(t)$, $\Pa \in \parset$ with $\Pa$ unbounded, there is a $\tilde{z}_t$ that belongs to the boundary of $\Pb(t)$, $\tilde{z}_t \in \boundary\Pb(t)$, such that $\norm{z_t - \tilde{z}_t} < \epsilon$. 
\end{assum}

\bigskip
Lemma~\ref{prop:exact_reachable_touching_trajectories} shows that any state $z_t \in \boundary \parinter(t)$ (with $t \in \Idef{\parset}$ given) is the terminal state of a touching trajectory $\Zas$ of $\parinter$ with initial state $\Zas(0) \in \boundary \parinter(0)$.
\begin{lem}
\label{prop:exact_reachable_touching_trajectories}
If Assumptions~\ref{hyp:boundary_reached} and \ref{hyp:dvg_paraboloid_regular} hold, any state $z_t \in \boundary \parinter(t)$ has a past touching trajectory $\Zas$ of $\parinter$ s.t. $\Zas(t) = z_t$.
\end{lem}
\begin{pf}
\newcommand{\lP}{\lambda}
\newcommand{\LP}{\overline{\lambda}}
Let $z_t \notin \reach(\parabol_0,t)$ such that for any $\epsilon>0$, there is a touching trajectory $\Zb$ of $\Pb$,  $\Pb$ finite, with  $\Zb(t) \in \reach(\parabol_0,t)$ and $\norm{\Zb(t) - z_t} < \epsilon$. For $\xqa(t) > 0$, we can define the optimal trajectory $\Zas$ with $\Zas(t) = \Za_t$. For any $\tau \in [0,t]$, $\Pa(\tau)$ is finite.
Then, Property~\ref{prop:case_psi} can be used over $[0,\tau]$.
Therefore, if $z_t \in \boundary \Pa$ such that $\Pa$ diverges at $t$, it holds:
$$
z_t \in \boundary \parinter(t) \Leftrightarrow z_t \in \reach(\parabol_0,t)
$$

For states not belonging to a diverging time-varying paraboloid, the property is a direct consequence of Assumption~\ref{hyp:boundary_reached}, Property~\ref{prop:case_psi}.
\end{pf}
Lemma~\ref{prop:exact_reachable_touching_trajectories} shows that any point on the boundary belongs to the reachable set since, for any given terminal state, we found a past trajectory (the touching trajectory) that satisfies the constraint \eqref{eq:cons} and with initial condition in the set of initial states.

\subsection{Exact Reachable Set}\label{ssec:exact_reachable_set}
We now state the main result of the paper:
\begin{thm}[Exact reachability, solution of Problem~\ref{pb:exact_reach_set}]
\label{thm:exact_reachable_set}
When Assumptions~\ref{hyp:boundary_reached} and \ref{hyp:dvg_paraboloid_regular} hold, the reachable set $\reach(\parabol_0,t)$ of system $\sys(\parabol_0,t)$ (defined in Section~\ref{sec:problem_statement}) is equal to the set $\parinter$ defined in \eqref{def:parinter}, namely
$$
\parinter(t) = \reach(\parabol_0,t)
$$
for all $t \geq 0$.
\end{thm}
\begin{pf}
Theorem~\ref{thm:overapproximation} states that $\reach(\parabol_0,t) \subseteq \parinter(t)$. By Property~\ref{prop:consume_energy}, for $z_t \in \parinter(t)$, we can construct a trajectory $\Za$ such that $\Za(t) = z_t$, $\Za(t^-) = z_t^* \in \boundary \parinter(t)$ (Property~\ref{prop:consume_energy}). Since $z_t^* \in \boundary \parinter(t)$, using Lemma~\ref{prop:exact_reachable_touching_trajectories}, there exists a trajectory $\Za$ such that $\Za(t^-) = z_t^*$ and $\Za$ is a touching trajectory of $\parinter$ on $\left[0,t\right[$. Since $\Za$ is a touching trajectory of $\parinter$, $\Za(0) \in \boundary \parinter(0)$ with $\parinter(0) = \parabol_0 = \reach(0)$. By Property~\ref{prop:reject}, the trajectory $\Za$ is valid (i.e. satisfies the energy constraint \eqref{eq:state_cons}) $z_t \in \reach(\parabol_0,t)$.
\end{pf}

\section{Implementation}
\label{sec:implementation}
\newcommand{\iparset}{\widetilde{\Pi}}%
\newcommand{\iscaleset}{\tilde{\Gamma}}%
\newcommand{\iparinter}{\widetilde{\Pi}^\cap}%
In this part, we discuss the practical implementation of the reachable sets overapproximation using Theorem~\ref{thm:exact_reachable_set}. To do so, we compute a subset $\iparset$ of $\parset$:
\begin{equation}
\label{eq:def_iparset}
\iparset \subseteq \parset
\end{equation}
$\iparset$ corresponds to the time-varying paraboloid set generated by a finite subset of scaling functions and initial scaling. Then, the intersection of each time-varying paraboloid evaluated at a given $t>0$ is an overapproximation of the reachable set $\reach(\parabol_0,t)$.
Finally, the DRE numerical integration is detailed for the case of non-negative solutions to the DRE.
We propose an algorithm (Algorithm~\ref{algo:approx}) that computes $\iparset$, its implementation in Matlab is available online \cite{mycode}.

\paragraph*{Subset of scaling functions and initial scaling factors:}
\newcommand{\Tcons}{T_c}
In this work, we choose to consider discrete scalings for the time-varying paraboloids.
The scalings are applied at each $k \Tcons$, for $\Tcons >0$ given and $k$ in $\N$.
A scaling is then described by a sequence of scaling factors $\{\lambda_k\}_{k\in \N}$, $\lambda_k \geq 1$, $k \in \N$. The scaling functions are not used: $\Ga(\cdot) = 0$. 

In the ideal case, the scaling (function and factors) would be chosen such that the following property is verified:
\begin{equation}
\exists \epsilon >0, \forall \tau \in [t,t+\epsilon], \xqas(\tau) \geq 0
\end{equation}
where $(\xas,\xqas)$ corresponds to the touching trajectory associated with the scaling function $\Ga$ and scaling factor $\ga$ and such that $(\xas(t),\xqas(t)) = (x,x_q)$.
In practice, since there might be an infinite number of states $(x,x_q)$ verifying $\dot{x}_q\geq0$, only a finite number of states are checked. These states are chosen as projections of a given point in given directions over $\Zps \cap \boundary \iparinter$.
These points are then used to evaluate a range of scaling factors $\ga$ to enforce $\dxqas(k \Tcons) \geq 0$. $\Ga$ is not used.

\paragraph*{Paraboloid numerical integration:}
Let two paraboloids $\Pa = \IVP(\parabol_0,1,\Ga)$, $\Pb = \IVP(\parabol_0,1,\Gb)$. If $\Ga(.) = \Gb(.)$ over an interval $[0,t_i]$, $t_i>0$, then $\Pa(.) = \Pb(.)$ over $[0,t_i]$.
Let $t_i \geq 0$ corresponds to the maximal time instant where there is $\Pb \in \parset$ such that $\Pb|_{[0,t_i]} = \Pa|_{[0,t_i]}$ (i.e. such that the restriction of $\Pb$ on $[0,t_i]$ is equal to the one of $\Pa$ on the same interval).
And let $t_f \geq 0$ corresponds either to the integration horizon $T>0$, or to the maximal of the interval of definition of $\Pa$.
For implementation purposes, each time-varying paraboloid is defined over the interval $[t_i,t_f] \subseteq [0,T]$.

Since $\Ga(.) = 0$ over $\left]k \Tcons,(k+1)\Tcons\right[$, for any $\tau \in [0,\Tcons]$, $k \in \N$. Over $\left]k \Tcons,(k+1)\Tcons\right[$, the IVP \eqref{eq:eq_diff_Efcg_gamma} is then equivalent to the "unscaled" IVP (\ref{eq:eq_diff_E},\ref{eq:eq_diff_fc},\ref{eq:eq_diff_g}) for $(\E_k,\fc_k,\r_k)$ over each time interval $\left[k \Tcons, (k+1) \Tcons\right]$ 
The solution to $\Pa = \IVP(\parabol_0,1,\Ga)$ is then described by parameters
$(\E,\fc,\g)$ with
$$(\E(t),\fc(t),\r(t)) = ( \E_k(t), \fc_k(t),  \r_k(t))$$
for each $t \in [k \Tcons, (k+1) \Tcons]$, $k \in \N$.

\paragraph*{Cardinal limitation of $\iparset$:}
In order to have a tractable integration of the reachable set computation, we limit the cardinality of $\iparset$ in the following way:
\newcommand{\Nnew}{N_{new}}
\newcommand{\Nparab}{N_{P}}
\begin{itemize}
\item at each time step $k \Tcons$, we consider only $\Nnew$ scaled paraboloids of highest scaling factor;
\item $\iparset$ below $\Nparab$, oldest time-varying paraboloids are dismissed in benefit of more recent ones;
\end{itemize}
$\Nnew$ and $\Nparab$ are user-defined parameters.
Choosing the paraboloids with this heuristic showed good results in practice. These rules try to only consider elements of $\parset$ that are more stable. Since for 2 solutions $\Ea$ and $\Eb$ of \eqref{eq:eq_diff_E} respectively defined over $[0,\Ta]$ and $[0,\Tb]$ where $\Ta,\Tb \in \R \cup \{\infty\}$, if $\Ea(0) \mleq \Eb(0)$, then $\Ea(t) \mleq \Eb(t)$ for $t$ in the interval of definition of $\Ea$ and $\Eb$, we have $\Ta \geq \Tb$ (these property follow directly by writing the corresponding value function of the basic LQR optimization problem). Therefore, for a time-varying paraboloid that is positive definite at $t>0$, its scaled time-varying paraboloid at $t$ will be defined for a longer time horizon.

\paragraph*{DRE numerical integration:} DRE integration is subject to numerical instability. A direct integration of the DRE~\eqref{eq:eq_diff_E} does not produce good results in practice (see \cite{kenney1985numerical}). Experiments presented in this works make use of the Chandrasekhar method \cite{Chandrasekhar}. This method integrates the Ordinary Differential Equation (ODE) \eqref{eq:eq_diff_E} $\E$ using an intermediate ODE over the time-dependent matrix $\sig{L}$ in $\Ld(\Rp,\R^{n \times n})$:
$$
\begin{aligned}
\dE(t) &= \sig{L}(t) \sig{L}(t)^\T\\
\dot{\sig{L}}(t) &= (\E(t) B \Mwi \Bt-\At-\Mxw \Mwi \Bt) \sig{L}(t)
\end{aligned}
$$
with 
$$
\begin{aligned}
\E(0) &= \E_0\\
\sig{L}(0) \sig{L}(0)^\T &= \dE_0
\end{aligned}
$$
where $\dE_0 = \dE(0)$ given by \eqref{eq:ode}. Then $\E$ is a solution to \eqref{eq:eq_diff_E}.

Since $\sig{L}(t) \sig{L}(t)^\T \mgeq 0$, this method is only applicable to strictly increasing solutions of the DRE. As seen in Example~\ref{ex:1dsys}, the solutions to ODE~\eqref{eq:eq_diff_E} are not strictly increasing over the time horizon, even for a positive definite initial condition. Therefore, the Chandrasekhar method cannot be used directly. We instead use the following approach, let $\sig{L},\sig{K} \in \Ld(\Rp,\R^{n \times n})$ such as:
$$
\begin{aligned}
\dE(t) &= \sig{L}(t) \sig{L}(t)^\T - \sig{K}(t) \sig{K}(t)^\T\\
\dot{\sig{L}}(t) &= (\E(t) B \Mwi \Bt-\At-\Mxw \Mwi \Bt) \sig{L}(t)\\
\dot{\sig{K}}(t) &= (\E(t) B \Mwi \Bt-\At-\Mxw \Mwi \Bt) \sig{K}(t)\\
\end{aligned}
$$
with
$$
\begin{aligned}
\sig{L}(0) \sig{L}(0)^\T &=  \dE_0^+\\
\sig{K}(0) \sig{K}(0)^\T &= -\dE_0^-
\end{aligned}
$$
where $\dE_0 = \dE(0) = \dE_0^+ + \dE_0^-$ given by \eqref{eq:eq_diff_E}, with $\dE_0^+ \mgeq 0$ and $\dE_0^- \mleq 0$  
The increasing and decreasing parts of $\E$ are respectively represented by the terms $\sig{L}$ and $\sig{K}$. Our Chandrasekhar inspired method performs better since the square root term $\sig{L}$ and $\sig{K}$ are much smaller than $\E$ and produces less numerical errors.

For $\fc$ and $\r$, integration of the ODE as given in \eqref{eq:eq_diff_fc} and \eqref{eq:eq_diff_g} is used.


Algorithm~\ref{algo:approx} summarizes the computation of $\iparset$. An implementation on Matlab is available online \cite{mycode}.

\newcommand{\textvar}[1]{\mathtt{#1}}
\IncMargin{1em}
\begin{algorithm}
\SetAlgoLined
\SetAlgoNoEnd
\DontPrintSemicolon
\Indm
\SetKwInOut{Input}{input}
\Input{\\
\hspace{-1cm}
\begin{tabular}[t]{p{7cm}}
A paraboloid $\parabol_0 \in \parabolSet$ of initial set of states\\
An horizon of simulation $T>0$\\
Sample time $\Tcons>0$ of constraint addition\\
Searching directions $\textvar{Search\_Dir} \subset \R^n$ to add constraints\\
$\textvar{\Nnew}$: maximal number of new scaled paraboloid to add\\
$\textvar{\Nparab}$: maximal cardinal of $\iparset$\\
\end{tabular}}
\KwResult{a set of overapproximating time-varying paraboloids $\iparset$}
\Indp
 $\iparset$ = $\{ \IVP(\parabol_0,1,0) \}$\;
 $t$ = 0\;
 $\textvar{Sim\_Parab} = \{(\parabol_0,0)\}$\;
 \While{$t < T$}{
 \BlankLine

 	\tcc{Find the new time-varying paraboloids to consider}
  $\textvar{New\_Parab} = \{\}$\;
  \For{$n \in \textvar{Search\_Dir}$}{
	   project $x_c$ on $\boundary \iparinter$ in the direction $n$\;
	   let $x^*$ be this projection and $P^* \in \iparset$ its corresponding touching paraboloid\;
       compute $\overline{\lambda}$ given ($x^*$,$P^*$)\;
	   
	  \For{$\lambda = 1+d\lambda, 1+2 d\lambda, \dots, \overline{\lambda}$}
	  {add $(P^*,\lambda)$ to $\textvar{New\_Parab}$\;}
	}  	
	  Sort $\textvar{New\_Parab}$ according to $\lambda$'s values\;
	  Keep $\textvar{\Nnew}$ elements of $\textvar{New\_Parab}$ with highest $\lambda$'s values
	  
	  \For{$(P^*,\lambda) \in \textvar{New\_Parab}$}{
		   add $(\lambda P^*(t),t)$ to $\textvar{Sim\_Parab}$\;
		}

		\If{$\card{\iparset}>\textvar{\Nparab}$}
		{
			Remove the $(\card{\iparset}-\textvar{\Nparab})$ oldest elements of $\textvar{Sim\_Parab}$\;
		}
}
	\BlankLine
\end{algorithm}

\begin{algorithm}%
\setcounter{AlgoLine}{17}%
\algocf@Vsline{%
 	\tcc{Simulate the paraboloid for $\Tcons$}
		\For{$(P_\tau,\tau) \in \textvar{Sim\_Parab}$}{
      Simulate $P(\cdot)$ over $[t, t+\Tcons]$ with $P(\tau)=  P_\tau$\;
      Add $P(\cdot)$ to $\iparset$\;
			\If{$P(\cdot)$ diverges}
			{
				Remove $(P,t)$ from $\textvar{Sim\_Parab}$\;
			} 
		}
		
		$t = t + \Tcons$\;
}
\caption{\label{algo:approx}Computation of $\iparset$ defined by \eqref{eq:def_iparset}, in Section~\ref{sec:implementation}, as the subset of $\parset$ defined by \eqref{def:parset}, in Section~\ref{sec:exact_reachable_set}.}
\end{algorithm}

\section{Examples}
\label{sec:examples}
Algorithm~\ref{algo:approx} deduced from Theorem~\ref{thm:overapprox_scaled} and \ref{thm:exact_reachable_set} is used to compute the overapproximation $\parinter$ defined in \eqref{eq:def_iparset} (subset of $\parset$ defined in \eqref{def:parset}) of reachable set $\reach(\parabol_0,t)$ of the system $\sys(\parabol_0,t)$ (described in Section~\ref{sec:problem_statement}), $t \geq 0$. Several examples are treated. With these examples, we provide some performance evaluations of our approach.

\subsection{Examples from COMPleib}
To evaluate the performance of our approach, we compute an overapproximation of the reachable set for several real-life systems from the COMPleib library \cite{leibfritz2006compleib}.
For each system, a stabilizing controller is generated for the generalized plant using the \texttt{h2syn} function of Matlab, then the system is reduced using a balanced truncation method to a given state space size.
The set of initial states is chosen such that the quadratic term belongs to the set of stable solutions to the associated Continuous Algebraic Riccati Equation.
The simulation are ran for an input $u(t) = \smat{1&\dots&1}^\T \exp(-t)$ for $t \in [0,2]$.
Each ODE is numerically integrated using the \texttt{ode113} solver in Matlab.
Finally, we run the simulation with one time-varying paraboloid and then multiple time-varying paraboloids. CPU time performances for a computer with an Intel i5 2.5GHz are presented in Table~\ref{table:cpu_perf}.

\begin{figure*}
\centering
\includegraphics[width=0.49\textwidth]{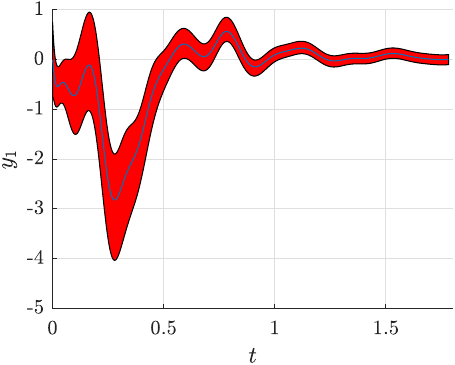}
\includegraphics[width=0.49\textwidth]{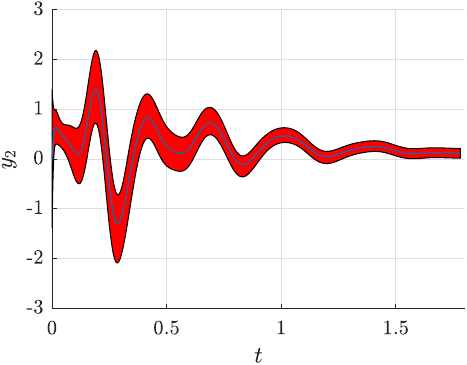}
\caption{\label{fig:ac10}Overapproximation of the output reachable set (projection of the reachable set $\R(t)$ through the observation map; the red area) of the AC10 example from the COMPleib library. Plain black line correspond to the unperturbed trajectory of the system.}
\end{figure*}

In Figure~\ref{fig:ac10}, we show several runs for the examples. Each paraboloid is overapproximated with a box, we show the intersection of these intervals.

Performance is mainly dependent on the number of paraboloids that we consider, and our ability to efficiently solve the DRE. 


\newcommand{\na}{\textit{n.a.}}
\begin{table}[h]
\centering
\newcommand{\tblsize}{2cm}
\begin{tabular}{p{1cm}|p{\tblsize}p{\tblsize}p{\tblsize}}
System size& Helicopter (HE7) & Aircraft (AC10) & Coupled Spring (CSE1) \\\hline
5   & 4.32 & 4.64  & 3.65  \\
10  & 5.12 & 5.96  & 3.86  \\
19  & 7.42 & 10.62 & 7.92  \\
30  & \na  & 28.85 & \na   \\
40  & \na  & 50.66 & \na   \\
49  & \na  & 88.00 & \na   \\
\end{tabular}
\caption{\label{table:cpu_perf}Computation times (in seconds) of the overapproximation for different systems sizes, using a unique time-varying paraboloid.
(When the original system's size is smaller than the required reduced system size, then the model reduction is \textit{not applicable} -\na-.)}
\end{table}

\begin{table}[h]
\centering
\newcommand{\tblsize}{2cm}%
\begin{tabular}{p{1cm}|p{\tblsize}p{\tblsize}p{\tblsize}}%
System size& Helicopter (HE7) & Aircraft (AC10) & Coupled Spring (CSE1) \\\hline
5   & 83.63 (66)  & 36.64 (13) & 213.88 (232) \\
10  & 89.55 (57)  & 25.77 (9)  & 261.32 (197) \\
19  & 167.53 (52) & 27.67 (4)  & 21.97 (4)    \\
30  & \na         & 113.72 (7) & \na          \\
40  & \na         & 117.60 (4) & \na          \\
\end{tabular}%
\caption{\label{table:cpu_perf}Computation times (in seconds) and number of paraboloids (in parenthesis) of the overapproximation for different systems sizes.
(When the original system's size is smaller than the required reduced system size, then the model reduction is \textit{not applicable} -\na-.)}
\end{table}

	\subsection{System Verification}
%
We study the stable IQC system
$\sys(\parabol_0,t)$, defined in \eqref{def:cons_sys}, at a given time $t$ in $[0,1]$, for a parabolic set of initial states 
$\parabol_0 = \parabol(\Ei_0,\fci_0,\ri_0)$, with $
\Ei_0 = \smat{a+b&a\\a&a+b}, \,
\fci_0 = \smat{0\\0}, \,
\ri_0 = 0.015, \,
a = 10^{-2}\textrm{ and }
b = 10^{-6}
$,
and for the following parameters 
$$A = -I, \,
B = I, \,
B_u = 0, \,
M = \smat{I&0&0\\0&1&0\\0&0&-2 I}
$$
where $I = \smat{1&0\\0&1}$, and with a zero input signal $u$.

\begin{figure*}[ht!]
\newcommand{\sizefig}{0.32\textwidth}%
\centering%
\begin{subfigure}[b]{\sizefig}
\includegraphics[width=\textwidth]{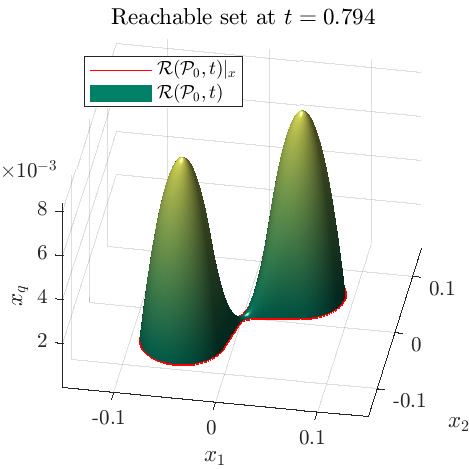}%
\caption{\label{fig:ex_stable_reach_set_3D}Reachable set}
\end{subfigure}
\begin{subfigure}[b]{\sizefig}
\includegraphics[width=\textwidth]{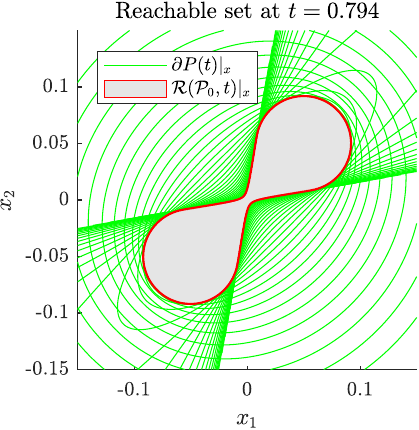}%
\caption{\label{fig:ex_stable_reach_set}Reachable set of the LTI system}
\end{subfigure}
\begin{subfigure}[b]{\sizefig}
\includegraphics[width=\textwidth]{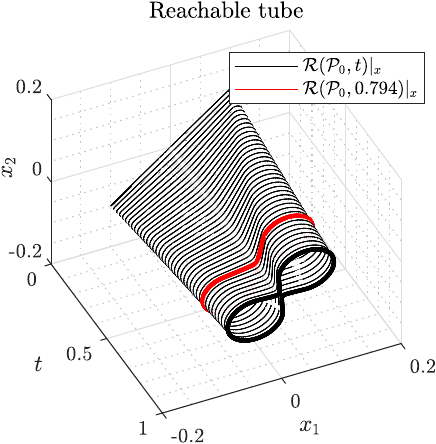}%
\caption{\label{fig:ex_stable_reach_tube}Reachable tube of the LTI system}
\end{subfigure}
\caption{The green surface in (a) is the reachable set $\reach(\parabol_0,t)$ at $t=0.794$ of $\sys(\parabol_0,t)$ computed using Theorem~\ref{thm:exact_reachable_set}. Its projection over the LTI state space $(x_1,x_2)$ (in solid red line) is shown in (b), each green line corresponds to one constraint $P \in \parset$ computed with Theorem~\ref{thm:overapproximation}. (c) is the reachable tube $t \rightarrow \reach(\parabol_0,t)$ of $\sys(\parabol_0,t)$ projected over the LTI state space $(x_1,x_2)$ for $t \in [0,1]$. The red section corresponds to the time $t = 0.794$.}
\end{figure*}

The reachable set  $\reach(\parabol_0,t)$ of $\sys(\parabol_0,t)$, defined in \eqref{eq:reach_set}, is computed using \eqref{def:exact_reach_set} and Theorem~\ref{thm:exact_reachable_set}, for $t \in [0,1]$.
Figures~\ref{fig:ex_stable_reach_set_3D} and \ref{fig:ex_stable_reach_set} show the reachable set $\reach(\parabol_0,t)$ set at time $t=0.794$ and its projection $\reach(\parabol_0,t)|_x$ over the LTI state space (i.e. projection over $(x_1,x_2)$ states).
In Figure~\ref{fig:ex_stable_reach_set}, the constraints boundaries $\boundary P(t)$ (for $P \in \parset$, $\parset$ defined in Section~\ref{sec:exact_reachable_set}) are touching the reachable set $\reach(\parabol_0,t)$. The non-convexity of $\reach(\parabol_0,t)$ arises from the non-positive solutions to the DRE~\eqref{eq:eq_diff_E}. Figure~\ref{fig:ex_stable_reach_tube} represents the projection of the reachable tube $t \mapsto \reach(\parabol_0,t)$ projected over the LTI dimension $(x_1,x_2)$.

\section{Conclusion}
\label{sec:conclusion}


IQCs are 2-norm constraints (i.e. energetic constraints) between signals.
Classical models in the robust control community involve 2-norm constraint or
$\infty$-norm constraint (i.e. hard bounds between signal). $\infty$-norm
bounds over the signals have been treated (e.g. the ellipsoidal method).
However, complex systems can be usually described with many relationships. In
future works, the computation of the reachable set for systems with multiple
2-norm and $\infty$-norm relationships will be investigated.  Also, such
models will be used to define abstractions for nonlinear systems. The dynamic
of the system will be linearized and the non-linearity modeled as a bounded
(2-norm and/or $\infty$-norm) disturbance.

\bigskip
We showed that the reachable set can be described as an intersection of
uncountably many paraboloids.  In our implementation, a subset of these
time-varying paraboloids is computed to overapproximate the reachable set. Then,
we compute a minimal volume paraboloid that contains the intersection of all the
paraboloids.  The computation time of our method is directly dependent of the
number of time-varying paraboloids. Finding only one time-varying paraboloid
which minimizes its end volume would avoid integrating multiple time-varying
paraboloids. Solutions exists for this optimization problem.

The differential Riccati equation can be weakly solved using a basis of
polynomial solutions. Then Sum-Of-Square relaxation provides a suboptimal
overapproximating paraboloid. Previous works implementing this approach use
conservative overapproximations that do not fully incorporate the state
constraint. In future works, we will develop such approach with the results
presented in this paper.

Locally optimal solution of the optimization problem can be derived using the
maximum of Pontryagin principle. Such solution are already available for the
ellipsoidal method. The adaptation to the paraboloidal method will be the topic
of a future work.

\bigskip
Our integration scheme is not guaranteed and the paraboloids we compute are
subject to the error of the differential equation numerical integration.
Guaranteed integration scheme exists to overapproximate the reachable set of
linear time-invariant systems. The image of an ellipsoidal set through the
matrix exponential is overapproximated. This result is then used to
overapproximate the reachable set of a linear system. In future works, we
will develop a similar approach for the paraboloidal method.

In our implementation, the scaling functions and initial scaling factor are
chosen such that some touching trajectory validate the constraint in the
future. Other criteria could be derived such as studying the average behaviors
of the trajectories. Since most of the computational effort are linear in the
number of time-varying paraboloids that needs to be simulated, an efficient
choice of the scaling factor can lead to algorithms that demand less
computational resources.

In control applications where the system is described by a partial differential
equation, a linear approximation of the model can be derived by projecting the
state over a finite basis of function. The approximation is then described by an
ordinary differential equation that usually have a high number of states
(several order of magnitudes). In its current implementation, our method proved
to be efficient for systems of less than a hundred states. For higher system
dimension, the numerical integration of the differential Riccati equation might
cumbersome. In future works, numerical integration of sparse differential
Riccati equation could be used to treat such examples.

\bibliography{mybiblio}

\begin{thebibliography}{10}
\expandafter\ifx\csname url\endcsname\relax
  \def\url#1{\texttt{#1}}\fi
\expandafter\ifx\csname urlprefix\endcsname\relax\def\urlprefix{URL }\fi
\expandafter\ifx\csname href\endcsname\relax
  \def\href#1#2{#2} \def\path#1{#1}\fi

\bibitem{blanchini2008set}
F.~Blanchini, S.~Miani, Set-theoretic methods in control, Springer, Boston:
  Birkh\"auser, 2008.

\bibitem{jaulin2001applied}
L.~Jaulin, M.~Kieffer, O.~Didrit, E.~Walter, Applied interval analysis: with
  examples in parameter and state estimation, robust control and robotics,
  Vol.~1, Springer-Verlag London, London, 2001.

\bibitem{bayen2007aircraft}
A.~M. Bayen, I.~M. Mitchell, M.~K. Osihi, C.~J. Tomlin, Aircraft autolander
  safety analysis through optimal control-based reach set computation, Journal
  of Guidance, Control, and Dynamics 30~(1) (2007) 68--77.

\bibitem{megretski1997system}
A.~Megretski, A.~Rantzer, System analysis via integral quadratic constraints,
  IEEE Transactions on Automatic Control 42~(6) (1997) 819--830.

\bibitem{megretski2010kyp}
A.~Megretski, {KYP} lemma for non-strict inequalities and the associated
  minimax theorem, arXiv preprint arXiv:1008.2552.

\bibitem{helmersson1999iqc}
A.~Helmersson, An {IQC}-based stability criterion for systems with slowly
  varying parameters, in: International Federation of Automatic Control, Vol.
  32: 14th World Congress, Elsevier, 1999, pp. 3183--3188.

\bibitem{megretski1997integral}
A.~Megretski, Integral quadratic constraints for systems with rate limiters,
  Tech. Rep. LIDS-P-2407, Massachusetts Institute of Technology, Laboratory for
  Information and Decision Systems, Cambridge, Massachusetts (1997).

\bibitem{peaucelle2009integral}
D.~Peaucelle, L.~Baudouin, F.~Gouaisbaut, Integral quadratic separators for
  performance analysis, in: European Control Conference, Budapest, 2009.

\bibitem{ariba2017}
Y.~Ariba, F.~Gouaisbaut, A.~Seuret, D.~Peaucelle, Stability analysis of
  time-delay systems via bessel inequality: A quadratic separation approach,
  International Journal of Robust and Nonlinear Control 28~(5) (2017)
  1507--1527.

\bibitem{chernousko1999}
F.~L. Chernousko, What is ellipsoidal modelling and how to use it for control
  and state estimation?, in: I.~Elishakoff (Ed.), Whys and Hows in Uncertainty
  Modelling, Springer, Vienna, 1999, pp. 127--188.

\bibitem{kurzhanski2002ellipsoidal}
A.~B. Kurzhanski, P.~Varaiya, On ellipsoidal techniques for reachability
  analysis. part {I}: external approximations, Optimization Methods and
  Software 17~(2) (2002) 177--206.

\bibitem{kurzhanskiy2007ellipsoidal}
A.~A. Kurzhanskiy, P.~Varaiya, Ellipsoidal techniques for reachability analysis
  of discrete-time linear systems, IEEE Transactions on Automatic Control
  52~(1) (2007) 26--38.

\bibitem{Chandrasekhar}
D.~Lainiotis, {Generalized Chandrasekhar algorithms: Time-varying models}, IEEE
  Transactions on Automatic Control 21~(5) (1976) 728--732.
\newblock \href {http://dx.doi.org/10.1109/TAC.1976.1101323}
  {\path{doi:10.1109/TAC.1976.1101323}}.

\bibitem{mycode}
{IQCARUS} matlab code, \url{https://github.com/roussePaul/IQCARUS}.

\bibitem{lee1967foundations}
E.~B. Lee, L.~Markus, Foundations of optimal control theory, John Wiley \&
  Sons, New York, 1976.

\bibitem{gusev2017extremal}
M.~I. Gusev, I.~V. Zykov, On extremal properties of boundary points of
  reachable sets for a system with integrally constrained control, in:
  International Federation of Automatic Control, Vol. 50: 20th World Congress,
  Elsevier, Toulouse, France, 2017, pp. 4082--4087.

\bibitem{soravia2000viscosity}
P.~Soravia, Viscosity solutions and optimal control problems with integral
  constraints, Systems \& Control Letters 40~(5) (2000) 325--335.

\bibitem{prajna2004safety}
S.~Prajna, A.~Jadbabaie, Safety verification of hybrid systems using barrier
  certificates, in: Hybrid Systems: Computation and Control, Springer, 2004,
  pp. 477--492.

\bibitem{henrion2014convex}
D.~Henrion, M.~Korda, Convex computation of the region of attraction of
  polynomial control systems, IEEE Transactions on Automatic Control 59~(2)
  (2014) 297--312.

\bibitem{korda2016moment}
M.~Korda, Moment-sum-of-squares hierarchies for set approximation and optimal
  control, Ph.D. thesis, EPFL, Switzerland (2016).

\bibitem{savkin1995recursive}
A.~V. Savkin, I.~R. Petersen, Recursive state estimation for uncertain systems
  with an integral quadratic constraint, IEEE Transactions on Automatic Control
  40~(6) (1995) 1080--1083.

\bibitem{guseinov2009approximation}
K.~G. Guseinov, Approximation of the attainable sets of the nonlinear control
  systems with integral constraint on controls, Nonlinear Analysis: Theory,
  Methods \& Applications 71~(1-2) (2009) 622--645.

\bibitem{jonsson2002robustness}
U.~J{\"o}nsson, Robustness of trajectories with finite time extent, Automatica
  38~(9) (2002) 1485--1497.

\bibitem{seiler2019finite}
P.~Seiler, R.~M. Moore, C.~Meissen, M.~Arcak, A.~Packard, Finite horizon
  robustness analysis of {LTV} systems using integral quadratic constraints,
  Automatica 100 (2019) 135--143.

\bibitem{LTVTools}
P.~Seiler, J.~Buch, R.~M. Moore, C.~Meissen, M.~Arcak, A.~Packard, {LTVTools
  (Beta)}, {A MATLAB Toolbox for Linear Time-Varying System} (2017.
  \url{https://www.mathworks.com/matlab central/fileexchange/69563-ltvtools}).

\bibitem{rousse:hal-02049158}
P.~Rousse, P.-L. Garoche, D.~Henrion, Parabolic set simulation for reachability
  analysis of linear time invariant systems with integral quadratic constraint,
  in: 18th European Control Conference, 2019, pp. 4301--4306.

\bibitem{CarstenJoost2018}
C.~W. Scherer, J.~Veenman, Stability analysis by dynamic dissipation
  inequalities: On merging frequency-domain techniques with time-domain
  conditions, arXiv e-prints.

\bibitem{savkin1996model}
A.~V. Savkin, I.~R. Petersen, Model validation for robust control of uncertain
  systems with an integral quadratic constraint, Automatica 32~(4) (1996)
  603--606.

\bibitem{kuvcera1973review}
V.~Ku{\v{c}}era, A review of the matrix {R}iccati equation, Kybernetika 9~(1)
  (1973) 42--61.

\bibitem{schuricht2000ordinary}
F.~Schuricht, H.~von~der Mosel, Ordinary differential equations with measurable
  right-hand side and parameters in metric spaces, Universit{\"a}t Bonn, 2000.

\bibitem{kenney1985numerical}
C.~Kenney, R.~Leipnik, Numerical integration of the differential matrix
  {R}iccati equation, IEEE Transactions on Automatic Control 30~(10) (1985)
  962--970.

\bibitem{leibfritz2006compleib}
F.~Leibfritz, {COMPleib: COnstrained Matrix optimization Problem library}
  (2006).

\end{thebibliography}

\end{document}